%% file: ms.tex
\documentclass[12pt]{article}
\usepackage{graphicx} 
\usepackage{hyperref}
\usepackage{longtable}
\hypersetup{
    colorlinks=true,
    citecolor=black,
    linkcolor=black,
    urlcolor=black
}

\usepackage{amsmath}
\usepackage{amssymb}
\usepackage{bbm}
\usepackage{geometry}

\usepackage{algorithm}
\usepackage{algpseudocode}

\usepackage[english]{babel}

\usepackage{booktabs}

\usepackage{xcolor}

\usepackage[numbers,super,comma,sort&compress]{natbib}

\usepackage{authblk}

\usepackage{tikz}
\usetikzlibrary{matrix}
\usetikzlibrary{external}
\input{commands}

\author[1]{Johannes Wiebe}
\author[2]{In\^es Cec\'ilio}
\author[1]{Ruth Misener}
\affil[1]{Department of Computing, Imperial College London, London, UK}
\affil[2]{Schlumberger Cambridge Research, Cambridge, UK}

\title{Robust optimization for the pooling problem}

\begin{document}

\maketitle

\begin{abstract}
    The pooling problem has applications, e.g., in petrochemical refining,
    water networks, and supply chains and is widely studied in global
    optimization.
    To date, it has largely been treated deterministically, neglecting the influence
    of parametric uncertainty.
    This paper applies two robust optimization approaches, reformulation and
    cutting planes, to the non-linear, non-convex pooling problem.
    Most applications of robust optimization have been either convex or
    mixed-integer linear problems.
    We explore the suitability of robust optimization in the context of
    global optimization problems which are concave in the uncertain
    parameters by considering the pooling problem with uncertain inlet
    concentrations. We compare the computational efficiency of reformulation
    and cutting plane approaches for three commonly-used uncertainty
    set geometries on 14 pooling problem instances and demonstrate how
    accounting for uncertainty changes the optimal solution.
    \begin{figure}[htb]
        \centering
        \includegraphics{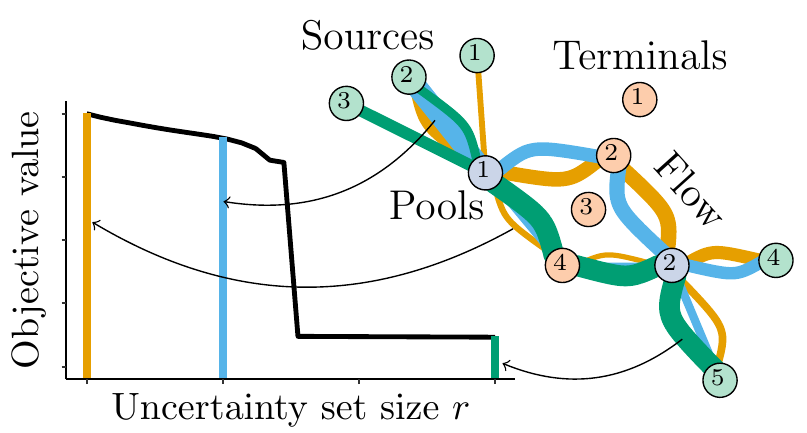}
    \end{figure}
\end{abstract}

\section{Introduction}

Robust optimization has become increasingly popular as a tool for planning and
scheduling of chemical and manufacturing processes under parametric uncertainty
\cite{Janak2005,Li2008,Bental2009,Gounaris2013,Zhang2015,Vujanic2016,Ning2017,Shang2018,Wiebe2018}.
The key assumption in robust optimization is that the uncertain parameters are
within some known, bounded uncertainty set.
Constraints containing the uncertain parameters are robustified by
requiring that they hold for all uncertainty set values.
Benefits of the robust approach are that
(i), {\corr for a rich class of problems and suitable uncertainty sets, considering
uncertainty induces only moderate increases in the number of variables
and constraints\citep{Bental2009,Bertsimas2010,MohajerinEsfahani2017a}} and
(ii) it does not require probabilistic information \citep{Bertsimas2011}
{\corr, even though such information may be used to construct uncertainty sets
with probabilistic guarantees \cite{Bental2009,Li2011b,Li2012b}}.

The robust optimization literature has traditionally studied
polynomially-solveable problems.
While starting with a focus on linear programming
problems (LPs)\cite{Soyster,Bental1999,Bertsimas2004}, robust optimization has
since evolved to encompass many non-linear
problem classes\cite{Bental1998,Bental2015,Marandi2017}.
In particular, {\corr a rich class of problems which are convex} in the
decision variables and concave in the uncertain parameters are often
tractable given a suitable uncertainty set\cite{Leyffer2018}.
In addition, the process systems engineering literature has used robust
optimization for many
applications in mixed-integer linear programming (MILP)\cite{Li2008,Li2011b}.

This work considers the connection of robust optimization and global
optimization in the context of the non-convex \emph{pooling problem}.
The pooling problem is relevant in oil and gas
refining, water systems, supply chains, and more\cite{misener-floudas:2009,Misener2011}.
It is non-linear and non-convex due to
bilinear mixing terms and is strongly $\mathcal{NP}$-hard\cite{Alfaki2013}.
The pooling problem is frequently used as a case
study in global optimization
\citep{Floudas1990,Meyer2006,Quesada1995,Tawarmalani2004,Wicaksono2008,
        Gounaris2009,Misener2010,Alfaki2013,Ceccon2016,Baltean2018}.
Parametric uncertainty in the pooling problem has been previously considered:
Li et al.\cite{Li2012a} apply robust optimization to demand
uncertainty in the linear product demand constraints while
Barton, Li, and co-workers develop a two-stage stochastic programming
approach for the extended pooling problem with uncertain inlet concentrations
and apply their methodology more generally to several applications
\cite{Li2011,Li2012,Li2013,Li2015,Yang2017,kannan-phd}.
A similar approach has also been proposed in the context of water network
synthesis\cite{Seong2014}.

{\corr
        While multi-stage approaches can lead to less
        conservative solutions because they can incorporate information which
        becomes available during operation, they tend to assume perfect knowledge
        of the uncertain parameters after the uncertainty has been revealed. In
        practice, real-time measurements may be able to reduce
        parametric uncertainty but not entirely eliminate it. In the context of
        pooling, the uncertain inlet concentrations may also change over time.
        Rigorously addressing these temporal changes would require a
        multi-stage, e.g., approximate dynamic programming, approach to
        managing the pooling
        network (due to accumulation in the pools). The single-stage robust
        approach offers an alternative to this which is especially
        relevant in applications, e.g., oil/natural gas networks or water
        networks, where some parametric uncertainty remains unresolved until
        the final blending step. This conservatism is relevant (i) when quality
        constraint violations may lead to unusable products or (ii) in health
        and safety-critical domains.
    }

The pooling problem is an interesting case study for global robust optimization,
because it is non-linear, non-convex in the decision variables but linear (and
therefore concave) in all potential uncertain parameters. Therefore many
robust optimization techiques, such as duality-based reformulations and robust
cutting plane approaches, are applicable. The resulting global optimization
problem is not significantly more difficult than the nominal pooling problem
without consideration of uncertainty.
While any of the problem's parameters could be considered uncertain, we focus on
the interesting case of uncertain component inlet concentrations, which occur
on the left hand side of the problems quality constraint in combination with
bilinear terms in the decision variables.
This work extends a previous conference paper\citep{Wiebe2019} which developed
robust reformulations for the pooling problem:
\begin{itemize}
    \item In addition to the duality-based robust reformulations from
        \citet{Wiebe2019}, we apply robust cutting
        planes to the non-convex pooling problem and develop a termination
        tolerance-based strategy for their integration with global
        optimization.
    \item We develop a simple but practical, safety factor inspired
        approach as a benchmark and as a computationally effective, conservative
        approximation to the robust problem.
    \item We assess and compare the suitability of these techniques in the context of global
        optimization and analyze their computational efficiency for several
        commonly used uncertainty set geometries (box, ellipsoidal, and
        polyhedral sets).
    \item We show how the optimal solution to the pooling problem is affected
        by varying degrees of uncertainty in the source component concentrations.
\end{itemize}

The first part of this paper reviews the different methods that have been
developed for the robust optimization problem and discusses the situations in
which each is applicable. It furthermore briefly discusses the integration of
global and robust optimization in the context of robust cutting planes. The
second part introduces the pooling problem, reviews its robust reformulations,
and develops the safety factor-inspired approach.
The final part compares the efficiency of the different methods for solving
the robust pooling problem and highlights the
effect of parametric uncertainty on the optimal flow configuration.

\section{Robust optimization}

Consider a generic robust optimization problem:
\begin{subequations}
    \begin{align}
        \min\limits_{\xvec\in \mathcal{X}} \; & f(\xvec)\\
        \text{s.t. } & g(\xvec, \xivec) \leq b && \forall \xivec \in
        \U,\label{cons:generic}
    \end{align}
    \label{prob:generic}
\end{subequations}
where $\xvec \in \R^n$ is the vector of decision variables, $\xivec \in \R^m$
the vector of uncertain parameters, $f(\cdot)$ the objective function,
$g(\cdot, \cdot)$ an uncertain constraint, $\U$ a non-empty, bounded set containing possible
realizations of $\xivec$, and $\X$ is the certain feasible region (which
does not depend on $\xivec$). Note that the objective does not depend on $\xivec$
and that we are only considering a single constraint. This formulation is
without loss of generality, because an uncertain objective can always be
transformed into an
uncertain constraint using an epigraph formulation and because robust optimization
{\corr usually considers each uncertain constraint separately (although
there are settings where constraints can be modeled
jointly\citep{Bental1998,Yuan2016})}.
In robust optimization, it is
often useful to rewrite Constraint~(\ref{cons:generic}) as an inner
maximization problem, leading to an equivalent bilevel
formulation\citep{Polak2012}:
\begin{subequations}
    \begin{align}
        \min\limits_{\xvec\in \mathcal{X}} \; & f(\xvec)\\
        \text{s.t. } & \max\limits_{\xivec \in \U} g(\xvec, \xivec) \leq b.
        \label{eq:inner-max}
    \end{align}
    \label{prob:bilevel}
\end{subequations}
The complexity of Problem~(\ref{prob:generic}) and the applicable solution
techniques largely depend on the
convexity/concavity of $\f$, $\g$, $\X$, and $\U$
with respect to $\xvec$ and $\xivec$ \citep{Leyffer2018}.
Four categories of problems exist:
\begin{enumerate}
    \item $\g$ and $\f$ are convex in $\xvec$ and concave in
        $\xivec$ and $\U$ and $\X$ are convex.
    \item $\g$ is concave in $\xivec$ and $\U$ is convex,
        but $\g$ or $\f$ is non-convex in $\xvec$ or $\X$ is non-convex.
    \item $\g$ and $\f$ are convex in $\xvec$ and $\X$ is convex,
        but $\g$ is non-concave in $\xivec$ or $\U$ is
        non-convex.
    \item $\g$ is non-concave in $\xivec$ or $\U$ is
        non-convex and $\g$ or $\f$ is non-convex in
        $\xvec$ or $\X$ is non-convex.
\end{enumerate}

\begin{figure}[htb]
    \centering
    \includegraphics{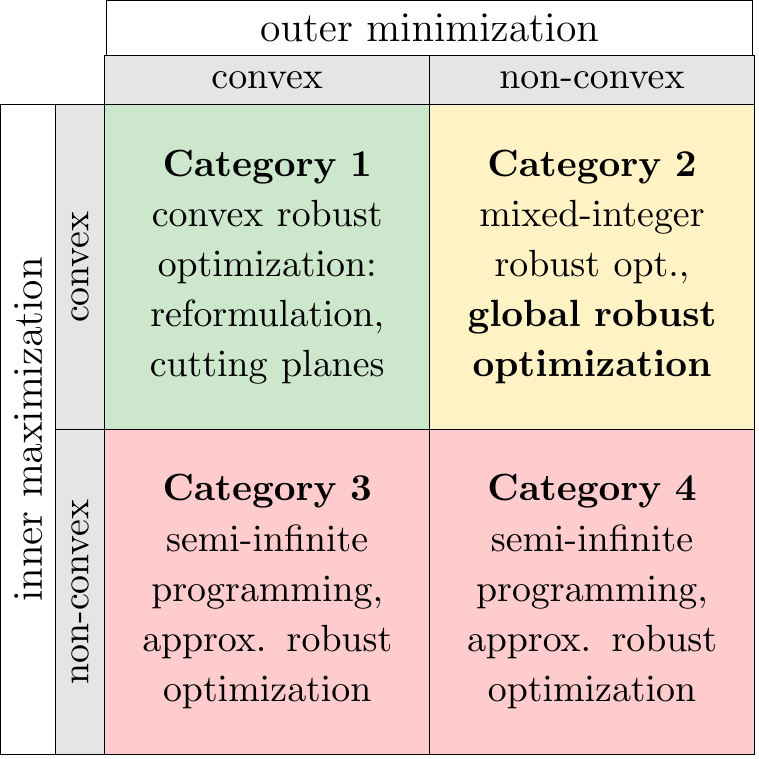}
    \caption{Complexity of Problem~(\ref{prob:generic}) for different convexity
    assumptions for the outer minimization and inner maximization in
    Problem~(\ref{prob:bilevel}). {\corr The inner maximization is non-convex when
$g(\cdot)$ is non-concave or $\U$ is non-convex.}}
    \label{fig:swot}
\end{figure}
Fig.~(\ref{fig:swot}) summarizes which methods are applicable to each category.
Most work in the robust optimization literature falls into
Category~1, also called convex robust optimization.
Category~1 includes linear programming problems with different types of
convex uncertainty sets, e.g., box sets \citep{Soyster}, ellipsoidal sets
\citep{Bental1999}, and polyhedral sets \citep{Bertsimas2004}.
It also includes non-linear convex problems, e.g., convex-quadratic problems
\citep{Bental1998,Marandi2017}, semi-definite programs \citep{Bental1998}, or
general convex constraints \citep{Ben-Tal2014}.
Both robust reformulation approaches and robust cutting planes have been
applied extensively Category~1 and solutions can often be
obtained in polynomial time.

Categories~3 and~4 are generally much harder, because the inner
maximization Problem~(\ref{eq:inner-max}) is a non-convex optimization problem.
Furthermore, even finding a feasible solution to
Problem~(\ref{prob:generic}) requires solving this non-convex problem to
global optimality.
Problems in these categories can 
generally only be solved using
semi-infinite programming techniques \citep{Tsoukalas2009,Mitsos2011,Stein2012},
which are limited to small scale problems, or approximate robust optimization
schemes \citep{Diehl2006,Houska2013,Yuan2017,Ho-nguyen2018}.

{\corr Category~2 has received limited attention in the robust optimization
and semi-infinite programming
communities\citep{Stein2003,Houska2012,Harwood2016}}.
A particular subclass, (non-convex)
MILP problems, has been studied extensively in the
process systems engineering literature, e.g.,
several robust planning and scheduling applications 
\citep{Li2008,Li2011b,Dias2016,Castro2018}.
We focus instead on the subclass of continuous but non-linear, non-convex
problems which are concave in the uncertain parameters, which includes the
pooling problem.
The same methods (both reformulation and cutting plane approaches)
which are used to solve problems in Category~1 are generally also applicable to
Category~2, but solutions can rarely be found in polynomial time.
Solving such problems effectively therefore requires a tighter integration between robust optimization and global optimization techniques.

\subsection{Reformulation approaches}
The robust reformulation approach was first proposed by Soyster\cite{Soyster}
for LPs with a box uncertainty set and has since been developed for many
classes of
problems\cite{Bental1998,Bental1999,Bertsimas2004,Ben-Tal2014,Marandi2017}.
These approaches generally utilize strong duality by replacing the inner
maximization in Problem~(\ref{prob:bilevel}) with its dual minimization problem:
\begin{subequations}
    \begin{align}
        \min\limits_{\xvec\in \mathcal{X}} \; & f(\xvec)\\
        \text{s.t. } & \min\limits_{\vec{\lambda} \in \Lambda} g^*(\xvec, \vec{\lambda}) \leq b,
    \end{align}
    \label{prob:dual}
\end{subequations}
where $g^*(\cdot, \cdot)$ is the objective of the dual problem, $\vec{\lambda}$
are the dual variables, and $\Lambda$ is the feasible space.
Problem~(\ref{prob:dual}) can clearly be written as a single level optimization
problem which can often be solved using off-the-shelf solver software:
\begin{subequations}
    \begin{align}
        \min\limits_{\xvec\in \mathcal{X},\vec{\lambda} \in \Lambda} \; & f(\xvec)\\
        \text{s.t. } & g^*(\xvec, \vec{\lambda}) \leq b.
    \end{align}
    \label{prob:rcr}
\end{subequations}
\subsection{Cutting plane approaches}

Robust cutting plane approaches outer-approximating the feasible space by replacing the
the uncertainty set $\U$ with a finite number of uncertainty
realizations $\Ucp_n = \left\{\xivec_1, \ldots,
\xivec_n\right\}$. A robustly feasible solution is found by iteratively
adding uncertainty realizations which violate the constraint.
The key idea is to alternately solve a master problem:
\begin{subequations}
\begin{align}
    \xvec_n = \argmax\limits_{\xvec \in \X} \; & f(\xvec)\\
    \text{s.t } & g(\xvec, \xivec) \leq b, \quad  \forall \xivec \in
    \Ucp_n,
\end{align}
    \label{prob:master}
\end{subequations}
and a separation problem:
\begin{align}
    \xivec_{n+1} = \argmax\limits_{\xivec \in \U} \; g(\xvec_n,
    \xivec),
    \label{prob:sep}
\end{align}
which determines whether the current solution $\xvec_n$ is robustly feasible
($\max\limits_{\xivec \in U} g(\xvec_n, \xivec) \leq 0$)
and, if it is not, generates a new cut $\xivec_{n+1}$.

Note that the master and separation problem are both convex optimization
problems when $\f$ and $\g$ are convex in $\xvec$
and concave in $\xivec$ and $\U$ and $\X$ are convex.
When $\f$ or $\g$ is non-convex in $\xvec$, e.g., in the pooling
problem, the master problem becomes a global
optimization problem which has to be solved at every iteration.
This could be disadvantageous because the global optimization solver may spend
a lot of time proving optimality of a solution which is not robustly feasible. One way to avoid this would be integrating the robust cutting planes
generation into the global optimization procedure.
The separation problem could be solved to
generate a new cut whenever an improved incumbent solution is found by the
global optimization solver, e.g., using lazy constraints.
These cuts have been proposed for MILP
\cite{Bertsimas2015}, but could also be useful in the global optimization
setting.

Unlike mature MILP solver software, however, global optimization solvers do not
routinely offer callbacks for adding lazy or incumbent constraints.
In light of this, we propose a simpler strategy for integration
of global optimization and robust cutting planes by adaptively changing the
termination tolerance for the master problem.
The general procedure is outlined in Algorithm~(\ref{algo:cp}).
\begin{algorithm}[t]
    \caption{Robust cutting planes}
    \begin{algorithmic}[1]
        \State $n \gets 0$, $\U_0 \gets \left\{\bar{\xivec}\right\}$,
        $\epsilon \gets \infty$
        \State $\delta \gets \delta_0$
        {\Comment \corr where $\delta$ is the master problem termination
        tolerance}
        \While{$\epsilon > 0 \text{ or } \delta > \delta^*$}
        \If{$\epsilon \leq 0$}
            \State $\delta \gets \delta \cdot \gamma$
            \Comment $\gamma < 1$
        \EndIf
        \State $
        \xvec_n \gets$ Eqn.~(\ref{prob:master})
        \Comment master problem
        \State $n \gets n + 1$
        \State $ \xivec_n\gets$ Eqn.~(\ref{prob:sep})
        \Comment separation problem
        \State $\Ucp_n \gets \Ucp_{n-1} \cup \left\{\xivec_n\right\}$
        \State $\epsilon \gets g(\xvec_n, \xivec_n)$
        \EndWhile
        \State \Return $\xvec_n$
    \end{algorithmic}
    \label{algo:cp}
\end{algorithm}
The master problem is first solved to a large tolerance $\delta_0$
Whenever a robustly feasible solution is found,
the tolerance $\delta$ is reduced by a factor $\lambda \leq 1$ until it reaches
the final tolerance $\delta^*$.

Robust cutting planes are similar to the concept of
\emph{backoff} in the control
literature\citep{Narraway1993,Bahri2006,Koller2018,Rafiei2018}.
The backoff approach generally also starts by
solving the nominal problem and then subsequently ``backing off'' to suboptimal
but more robust solutions which remain feasible in light of uncertainty.
While the robust optimization approach assumes that the uncertain parameter is
within a known uncertainty set, uncertainty in the backoff approach is usually
due to process dynamics.

\section{Robust pooling}
\label{sec:rob}
Consider the $q$-formulation of the standard pooling problem
\citep{haverly,bental}:
{\allowdisplaybreaks
\begin{subequations}
\begin{align}
        & \min\limits_{\q, \y, \z}
        &&\begin{array}{l} \sum\limits_{\substack{%
                            (i,l) \in T_X\\
                            (l,j) \in T_Y}} c_i \cdot \q \cdot \y
        - \sum\limits_{(l,j) \in T_Y} d_j \cdot y_{l,j}\\
        - \sum\limits_{(i,j) \in T_Z} (d_j - c_i)\cdot \z
        \end{array}
        && \;
        \label{eq:q-obj} \\
        &\begin{array}{l}
            \text{Feed}\\ \text{Availab.}
        \end{array}
        && \left\{
            A^L_i \leq
            \begin{array}{l}
            \sum\limits_{\substack{%
                            l:(i,l) \in T_X\\
                            (l,j) \in T_Y
                    }} q_{i,l} \cdot y_{l,j}  \\
            + \sum\limits_{j:(i,j) \in T_Z} \z
            \end{array}
            \leq A_i^U,
        \right.
        && \forall i \label{eq:q-avail}\\
        &\begin{array}{l}
            \text{Pool}\\ \text{Capacity}
        \end{array}
        && \left\{\sum\limits_{j:(l,j) \in T_Y} y_{l,j} \leq S_l, \right.
        && \forall l\\
        &\begin{array}{l}
            \text{Product}\\ \text{Demand}
        \end{array}
        &&\left\{ D_j^L \leq
        \begin{array}{l}
            \sum\limits_{l:(l,j) \in T_Y} y_{l,j} \\
            + \sum\limits_{i:(i,j) \in T_Z} z_{i,j}
        \end{array} \leq D_j^U,
        \right.
        && \forall j\\
        & \begin{array}{l}\text{Simplex}\end{array}
        && \left\{ \sum\limits_{i:(i,l) \in T_X} \q = 1,\right.
        && \forall l \label{eq:q-simplex}\\
        & \begin{array}{l} \text{Product} \\ \text{Quality} \end{array}
        &&
        \left\{ \begin{array}{l}
                    \sum\limits_{\substack{%
                                    l:(l,j) \forall T_Y\\
                                    i:(i,l) \forall T_X}} \nompar \cdot \q
                                    \cdot \y\\
                    + \sum\limits_{i:(i,j) \in T_Z} \nompar \cdot \z\\
                \end{array}
                \begin{cases}
                    \geq \PL \cdot \vj\\
                    \leq \PU \cdot \vj
                \end{cases}
            \right.
            && \forall j,k. \label{eq:q-quality}
\end{align} \label{prob:q-form}
\end{subequations}
}
The decision variables are the fraction of flow $q_{i,l}$ from source $i$ to pool $l$,
the flow $y_{l,j}$ from pool $l$ to terminal $j$, and the
direct flow $z_{i,j}$ from source $i$ to terminal $j$.
The sets $T_X, T_Y, T_Z$ describe the connections between sources,
pools and terminals.
$\AL/\AU$, $\DL/\DU$ and $\PL/\PU$ are the lower and upper limits for feed
availability, product demand, and product quality respectively.
The pool capacity is limited by $S_l$, the
inlet concentration of component $k$ at source $i$ is $\nompar$, product $j$ is
sold at price $d_j$ and the cost of material from source $i$ is $c_i$.
{\corr The objective minimizes negative profit.}
For notational brevity we introduce the total flow $\vj$ at terminal $j$ in
Constraint~(\ref{eq:q-quality}):
\begin{equation*}
    \vj = \sum\limits_{l:(l,j) \in T_Y} \y + \sum\limits_{i:(i,j) \in T_Z} \z \quad \forall j.
\end{equation*}
While non-linearities (bilinear mixing terms) in the decision variables occur in the
objective function, feed availability, and product quality constraints, the
problem is linear in all (potentially uncertain) parameters.
The problem therefore belongs to Category~2 outlined above, independent of
which parameter is considered to be uncertain.
Furthermore, unlike the original $p$-formulation of the pooling problem,
there are no potentially uncertain equality constraints. In fact, the only equality
constraint is the simplex which does not contain parameters. This makes the
$q$-formulation particularly suitable for robust optimization, because
uncertain equality constraints generally either have to be eliminated or treated using
adjustable robust optimization \cite{Gorissen2015,Lappas2016}.

While we could consider any or all of the problems parameters to be uncertain,
we focus on robustifying the product quality Constraints~(\ref{eq:q-quality})
considering uncertain input component concentrations $\uncpar$:
\begin{equation}
            \sum\limits_{\substack{%
                            l:(l,j) \forall T_Y\\
                            i:(i,l) \forall T_X}} \uncpar \cdot \q
                            \cdot \y
            + \sum\limits_{i:(i,j) \in T_Z} \uncpar \cdot \z
        \begin{cases}
            \geq \PL \cdot \vj\\
            \leq \PU \cdot \vj
        \end{cases}
        \forall \uncpar \in \Up, \; \forall j,k. \label{eq:rob-cons}
\end{equation}
This constraint is particularly interesting, because the uncertain parameters
occur on the left hand side and in combination with the bilinear terms
$\q\cdot\y$.

\subsection{Uncertainty set}
We consider a commonly used class of uncertainty sets $\Up$ defined by the
$p$-norm\cite{Li2011b}:
\begin{equation}
    \Up = \left\{\uncpar = \nompar + \devpar \cdot \xi_{i,k} \; | \;
    \|\xivec_k\|_{p} \leq r \right\},
    \label{eq:Up}
\end{equation}
where $\xivec_k = [\ldots, \xi_{i,k}, \ldots]^T$, \nompar is the nominal value of
\uncpar, \devpar its maximum deviation, and $r$ a size parameter which
determines how much of the uncertainty is considered (For $r = 0$ the
constraint is equivalent to the nominal constraint and for $r = 1$ all possible
realizations of $\uncpar$ are in the uncertainty set. The set $\U_p$ corresponds
to a box uncertainty set for $p=\infty$, an ellipsoid for $p=2$, and a
polyhedron for $p=1$.
All of these uncertainty set geometries have been discussed extensively
in the robust optimization literature~\cite{Li2011b}.

{\corr
While these uncertainty sets can capture different types of uncertainty
scenarios, they do not utilize any known correlation between the random
variables, e.g., the geographic structure of the pooling network.
For example, inlet concentrations at geographically neighboring sources may be
correlated more strongly than geographically distant sources.
To leverage such dependence structure, we constructe ellipsoidal uncertainty
sets based on distance between sources
using kernel functions. Kernel functions, e.g., the squared exponential kernel:
\begin{equation}
    k(\boldsymbol{\ell}, \boldsymbol{\ell}') = \sigma^2 \exp\left(\frac{\|\boldsymbol{\ell} - \boldsymbol{\ell}'\|_2^2}{2l^2}\right)
    \label{eq:se-kernel}
\end{equation}
where $\sigma^2 > 0$ is the signal variance and $l > 0$ is the length scale, can be used to model correlation between
two geographic locations $\boldsymbol{\ell}$ and $\boldsymbol{\ell}'$\citep{Chiles2011}.
We construct a
covariance matrix $\Sigma$ with elements $\sigma_{i, i'}$ describing the
correlation between each pair of
network sources $i$ and $i'$:
\begin{equation*}
    \sigma_{i,i'} = k(\boldsymbol{\ell}_i, \boldsymbol{\ell}_{i'}),
\end{equation*} where $\ell_i$ is the location of source $i$.
Note that this approach could be combined
with Gaussian process regression (Kriging) to construct $\Sigma$
 based on available concentration data from different
locations\citep{Chiles2011}.
Using $\Sigma$, Eqn.~(\ref{eq:u-corr}) constructs
an ellipsoidal uncertainty set exploiting the correlation between sources:
\begin{equation}
    \left\{\uncpar = \nompar + \devpar \cdot \xi_{i,k} \; | \; \xivec_k\T
        \Sigma^{-1} \xivec_k  \leq r^2\right\}.
        \label{eq:u-corr}
\end{equation}
We note that, if $\sigma = 1$ in Eqn.~(\ref{eq:se-kernel}) and the sources
are very far from each other, i.e., there is no correlation, $\Sigma$ becomes
the identity matrix and the uncertainty set in Eqn.~(\ref{eq:u-corr}) is
equivalent to $\U_2$ (Eqn.~\ref{eq:Up}).
}

\subsection{Reformulation}

Because Constraint~(\ref{eq:rob-cons}) is linear in the uncertain parameters
$\uncpar$ and $\Up$ is convex, duality-based robust reformulation techniques are
applicable.
To this end, we introduce the indicator function:
\begin{equation*}
    \mathbbm{1}(s \in S) =
    \begin{cases}
        1 & s \in S\\
        0 & s \notin S,
    \end{cases}
\end{equation*}
allowing us to rewrite Constraint~(\ref{eq:rob-cons}) in standard linear form:
\begin{equation}
    \begin{array}{ll}
        \sum\limits_{i}{\uncpar \cdot \x}
        \begin{cases}
            \geq \PL \cdot \vj\\
            \leq \PU \cdot \vj,
        \end{cases}
        & \forall \uncpar \in \Up, \forall j,k,
    \end{array}
    \label{eq:rob-linear}
\end{equation}
where $\x$ is the total flow from source $i$ to terminal $j$:
\begin{equation}
    \x = \sum_{l:(l,j) \in T_Y}
    \left[%
        \mathbbm{1}\left((i,l) \in T_X\right) \cdot y_{l,j} \cdot q_{i,l}
    \right]
    + \mathbbm{1}\left((i,j) \in T_Z\right) \cdot z_{i,j}.
    \label{eq:x}
\end{equation}

Using robust optimization results \cite{Li2011b}, the semi-infinite
Constraint~(\ref{eq:rob-linear}) can be rewritten as an equivalent
deterministic constraint:
\begin{align}
    & \sum_{i} \nompar \cdot \x
        \begin{cases}
            \leq \PU \cdot \vj - r \cdot \hp\\
            \geq \PL \cdot \vj + r \cdot \hp
        \end{cases}
    && \forall j, k
        \label{eq:rob-reform}
\end{align}
where $\hp$  depends on the selected uncertainty set geometry $p$
as indicated in Table~(\ref{tab:padding}).
\begin{table}[t]
    \centering
    \begin{tabular}{lll}
        \toprule
        Uncertainty set geometry & $p$ & $\hp$\\
        \midrule
        Box         & $\infty$  & $\sum\limits_i \devpar \cdot \x$ \\
        Ellipsoid   & $2$       & $\sqrt{\sum\limits_i \devpar^2 \cdot \x^2}$ \\
        \corr
        Ellipsoid (correlation)&\corr $2$       &\corr $\sqrt{
                                                \sum\limits_i
                                                \sum\limits_{i'}
                                                \sigma_{i, i'} \cdot \devpar
                                                \cdot \hat{C}_{i', k} \cdot \x
                                                \cdot x_{i', j}}$\\
        Polyhedron  & $1$       & $\max\limits_{i} \devpar \x$ \\
        \bottomrule
    \end{tabular}
    \caption{Robust reformulations for box, ellipsoidal, and polyhedral
        uncertainty sets.}
    \label{tab:padding}
\end{table}
Note that no absolute values are required because the flow $\x$ is always
positive.

For the box set substituting Eqn.~\ref{eq:x} for $\x$ leads to:
\begin{align}
    \begin{array}{l}
        \sum\limits_{\substack{%
            l:(l,j) \in T_Y\\
            i:(i,l) \in T_X
        }}%
        \nompar \cdot \q \cdot \y\\
        + \sum\limits_{i:(i,j) \in T_Z} \nompar \cdot z_{i,j}
    \end{array}
    \begin{cases}
        \leq
        \PU \cdot \vj -
        r \left[ \sum\limits_{\substack{%
        l:(l,j) \in T_Y\\
        i:(i,l) \in T_X
    }}%
    \devpar \cdot \q \cdot \y
+ \sum\limits_{i:(i,j) \in T_Z} \devpar \cdot z_{i,j}\right]\\
    \geq \PL \cdot \vj  + r\left[
    \sum\limits_{\substack{%
        l:(l,j) \in T_Y\\
        i:(i,l) \in T_X
    }}%
    \devpar \cdot \q \cdot \y
+ \sum\limits_{i:(i,j) \in T_Z} \devpar \cdot z_{i,j} \right]
    \end{cases}
    \label{eq:box-reform}
\end{align}
which simplifies to:
    \begin{align*}
        \sum\limits_{\substack{%
            l:(l,j) \in T_Y\\
            i:(i,l) \in T_X
        }}%
        \left(\nompar + r\cdot\devpar\right)\cdot \q \cdot \y
        + \sum\limits_{i:(i,j) \in T_Z} \left(\nompar + r\cdot\devpar\right) \cdot z_{i,j}
            \leq
            \PU \cdot \vj \quad \forall j, k,
    \end{align*}
    and
    \begin{align*}
        \sum\limits_{\substack{%
            l:(l,j) \in T_Y\\
            i:(i,l) \in T_X
        }}%
        \left(\nompar - r\cdot\devpar\right)\cdot \q \cdot \y
        + \sum\limits_{i:(i,j) \in T_Z} \left(\nompar - r\cdot\devpar\right) \cdot z_{i,j}
            \geq
            \PL \cdot \vj \quad \forall j, k,
    \end{align*}
respectively. In the simple case of the box uncertainty set, the worst case
uncertainty scenario is clearly always $\nompar + r\cdot\devpar$ for the upper
quality constraint and $\nompar - r\cdot\devpar$ for the lower quality
constraint.
The reformulation of the ellipsoidal uncertainty set can potentially be improved by
squaring each side of the quality constraints:
\begin{align*}
    r^2 \sum\limits_i \devpar^2 \x^2
    \begin{cases}
        \leq (\PU \cdot \vj - \sum\limits_i \nompar \cdot \x)^2\\
        \leq (\sum\limits_i \nompar \cdot \x - \PL \cdot \vj)^2,
    \end{cases}
\end{align*}
This constraint eliminates the square root which can be numerically challenging, but is
only valid when the nominal Constraint~(\ref{eq:q-quality}) is also added to the
model, increasing the number of equations.
{\corr The same applies to the ellipsoidal uncertainty set with correlation.}
For the polyhedral set, Constraint~(\ref{eq:rob-cons}) can be written as:
\begin{subequations}
    \begin{align}
        \sum\limits_{\substack{%
            l:(l,j) \in T_Y\\
            i:(i,l) \in T_X
        }}%
        \nompar \cdot \q \cdot \y
        + \sum\limits_{i:(i,j) \in T_Z} \nompar \cdot z_{i,j}
        &
        \begin{cases}
            \leq \PU \cdot \vj - r \cdot \hp\\
            \geq \PL \cdot \vj  + r \cdot \hp
        \end{cases} && \forall j, k\\
        \hp  & \geq \devpar \cdot \x && \forall i, j, k
    \end{align}
\end{subequations}

While the box uncertainty set does not change the problem complexity
with respect to the nominal case, the ellipsoidal and polyhedral
uncertainty sets both add (potentially many) convex constraints to the model.
These constraints do not change the complexity class, but they may introduce practical difficulties.
While the ellipsoidal set adds only $\mathcal{O}(|i|\cdot|j|)$ constraints,
the polyhedral set adds $\mathcal{O}(|i| \cdot |j| \cdot |k|)$ constraints and
variables. An advantage of the
polyhedral set, however, is that the problem remains bilinear, while the
ellipsoidal set introduces higher order terms.

\subsection{Cutting planes}
The cutting plane approach outlined above is directly applicable to the pooling
problem. The master problem is:
\begin{subequations}
\begin{align}
    & \min\limits_{\q, \y, \z} && \text{Eqn.~(\ref{eq:q-obj})}\\
    & \text{s.t. } && \text{Eqns.~(\ref{eq:q-avail}-\ref{eq:q-simplex})}\\
    &&&               \text{Eqn.~(\ref{eq:x})} \\
    &
    && \sum\limits_i \uncpar \cdot \x
                \begin{cases}
                    \geq \PL \cdot \vj\\
                    \leq \PU \cdot \vj
                \end{cases}
    &&& \forall \uncpar \in \Ucp_n, \forall j, k
\end{align}
\end{subequations}
where $\Ucp_0 = \left\{\nompar\right\}$ contains only the nominal values.
The separation problem is:
\begin{equation}
    \max\limits_{j, k}
    \max\limits_{\uncpar \in \U_p}
    \max \left\{
        \sum\limits_i \uncpar \cdot \x - \PU \cdot \vj,
        \PL \cdot \vj -\sum\limits_i \uncpar \cdot \x 
    \right\},
    \label{prob:pool-sep}
\end{equation}
which can be solved by solving $2\cdot|J|\cdot|K|$ convex optimization
problems. For the uncertainty sets considered in this work, it can also be
reformulated as a MILP (for box and polyhedral sets) or a mixed-integer
quadratically constrained program (MIQCP) using a big-M reformulation which can
be solved effectively with existing MIP solvers.

We explore two different ways of adding cuts: adding a single cut in every
round, consisting of the $j^*$-$k^*$-quality constraint which is most
violated, and adding multiple cuts in every round, consisting of all
$2\cdot|J|\cdot|K|$ quality constraints for the worst case scenario $\uncpar^*$,
where ($j^*, k^*, \uncpar^* = \argmax$ Eqn.~(\ref{prob:pool-sep})).
{\corr
    Note that an intermediate approach which adds the first $n$ most
    violated cuts would also be possible.
    Furthermore, $\Ucp_0$ could also be initialized with other values. For
    example, for the box uncertainty set the maximum and minimum concentrations
    could be used. This would make the master problem essentially equivalent to
    the reformulation (Eqn.~\ref{eq:box-reform}) and the cutting plane algorithm
    would always converge after the first iteration.
}

{\corr
    While there are examples for which Algorithm~\ref{algo:cp} does not
    converge \citep{Hettich1993}, it often performs similar to the reformulation
    approach in practice \citep{Fischetti2012,Bertsimas2015}. Furthermore, the
    algorithm is guaranteed to converge for the robust pooling problem with the Table \ref{tab:padding} uncertainty sets.
    Since Eqn.~\ref{eq:inner-max} is linear in $\uncpar$, the separation
    problem generally has a finite number of solutions for the box and
    polyhedral sets, the vertices of the polytope. For the ellipsoidal set, M\'inguez and
    Casero--Alonso\citep{Minguez2019} show that convergence is
    guaranteed when the constraint obtained with the reformulation approach is
    convex by showing that the cutting plane approach is equivalent to solving
    the reformulation using outer-approximation. Since all of the
    non-convexities are handled by the global optimization solver and the
    robust reformulation (Eq.~\ref{eq:rob-reform}) is convex in $\x$, the M\'inguez and
    Casero--Alonso\citep{Minguez2019} result also applies to the pooling problem.
}

\subsection{Safety factors approach}
A very simple but practically relevant approach to design under uncertainty is
the concept of safety factors \citep{Dao1974,Clausen2006,Hansson2010}. In this
approach, the space of feasible designs is restricted by scaling bounds with a
fixed safety factor $s > 0$. In the case of the pooling problem with uncertain
inlet concentration, the upper (and lower) bound $\PU$ (and $\PL$) may be scaled to discard
solutions which may easily become infeasible, while $\uncpar$ is replaced by its
nominal value $\nompar$:
\begin{subequations}
\begin{align}
    & \min\limits_{\q, \y, \z} && \text{Eqn.~(\ref{eq:q-obj})}\\
    & \text{s.t. } && \text{Eqns.~(\ref{eq:q-avail}-\ref{eq:q-simplex})}\\
    &
    &&
            \sum\limits_{\substack{%
                            l:(l,j) \forall T_Y\\
                            i:(i,l) \forall T_X}} \nompar \cdot \q
                            \cdot \y
            + \sum\limits_{i:(i,j) \in T_Z} \nompar \cdot \z
        \begin{cases}
            \geq \left(\PL\cdot s\right) \cdot \vj\\
            \leq \left(\PU\cdot \frac{1}{s}\right) \cdot \vj.
        \end{cases}
\end{align}\label{prob:sf}
\end{subequations}
A benefit of safety factors is that they do not increase the complexity of the
problem at all.
A disadvantage is, however, that one generally has to rely on experience when
choosing a suitable value for $s$. Note that Problem~(\ref{prob:sf}) is
equivalent to a robust version of Problem~(\ref{prob:q-form}) with uncertain
bounds $\PU$ and $\PL$ and a box uncertainty set, e.g.:
\begin{equation*}
    \PUunc \in
    \left\{\PUunc \; |  \; \PU \cdot \frac{1}{s}
             \leq \PUunc
             \leq \PU \cdot s\right\}.
\end{equation*}
For $s = 1$, Problem~(\ref{prob:sf}) is equivalent to the nominal problem and
for $s \to \infty$ the optimal solution is 
$\y = \z = 0$ for all $(l,j)$, $(i,j)$, i.e., the safest solution is producing nothing.

To compare the safety factor approach with the robust optimization
approach, we find the smallest safety factor $\smin$ which gives a
feasible solution for the robust problem with a
given uncertainty set $\U_p$. This is outlined in Algorithm~(\ref{algo:sf}).
We employ a combination of bisection and secant search to find the safety
factor for which Eqn.~(\ref{prob:pool-sep}) is just feasible.
\begin{algorithm}[t]
    \caption{Optimal safety factor approach}
    \begin{algorithmic}[1]
        \State $s_0 \gets 1, {\corr s_1 \gets \bar{s}}$
        {\Comment \corr where $\bar{s}$ is a large-enough number, e.g., 100}
        \State $\delta \gets 10^{-6}$
        \State $(\y, \z, \q)_{0/1} \gets \argmin$ Eqn.~(\ref{prob:sf})
        \Comment safety factor problem with $s_0$ and $s_1$
        \State $\epsilon_{0/1} \gets$ Eqn.~(\ref{prob:pool-sep})
        \Comment separation problem with solution from Eqn.~(\ref{prob:sf})
        \While {$|s_1 - s_0|/s_0 \geq \delta$}
            \If{$\epsilon_1 \leq 0.01$}
                \State $s \gets \frac{s_0 + s_1}{2}$
                \Comment bisection method when close to feasible
            \Else
                \State $s \gets s_0 + \epsilon_0 \frac{s_0 - s_1}{\epsilon_1 -
                        \epsilon_0}$
                \Comment secant method otherwise
            \EndIf

            \State $(\y, \z, \q) \gets$ Eqn~(\ref{prob:sf})
            \Comment safety factor problem
            \State $\uncpar \gets \argmax$ Eqn~(\ref{prob:pool-sep})
            \Comment separation problem
            \If{$(\y, \z, \q)$ is a feasible solution to Eqn.~(\ref{eq:rob-cons})}
                \State $s_0 \gets s$, $\epsilon_0 \gets$
                Eqn.~(\ref{prob:pool-sep})
            \Else
                \State $s_1 \gets s$, $\epsilon_1 \gets$
                Eqn.~(\ref{prob:pool-sep})
                \State $\smin \gets s_1$
            \EndIf
        \EndWhile
        \State \Return $(\y, \z, \q)$
    \end{algorithmic}
    \label{algo:sf}
\end{algorithm}

\section{Results}
\label{sec:res}

The robust reformulation, robust cutting plane, and safety-factor approach
outlined above were applied to 14
literature pooling instances \citep{adhya,foulds,bental,haverly,rt2}.
The model was implemented in GAMS 26.1.0 and solved on an i7-6700 CPU with
$8\times3.4$GHz and 16GB RAM.
NLP problems were solved using ANTIGONE 1.1 \citep{antigone} while MILP
problems were solved using CPLEX 12.8.
Each instance was solved for 30 values of $r$ for the box, ellipsoidal,
and polyhedral uncertainty set.
A time limit of 1~hr, a relative termination tolerance of $10^{-6}$,
and $\devpar = \nompar$ were used. For the cutting plane approach, a maximum
of 200 cuts were added.

\begin{figure}[htb!]
    \centering
    \includegraphics{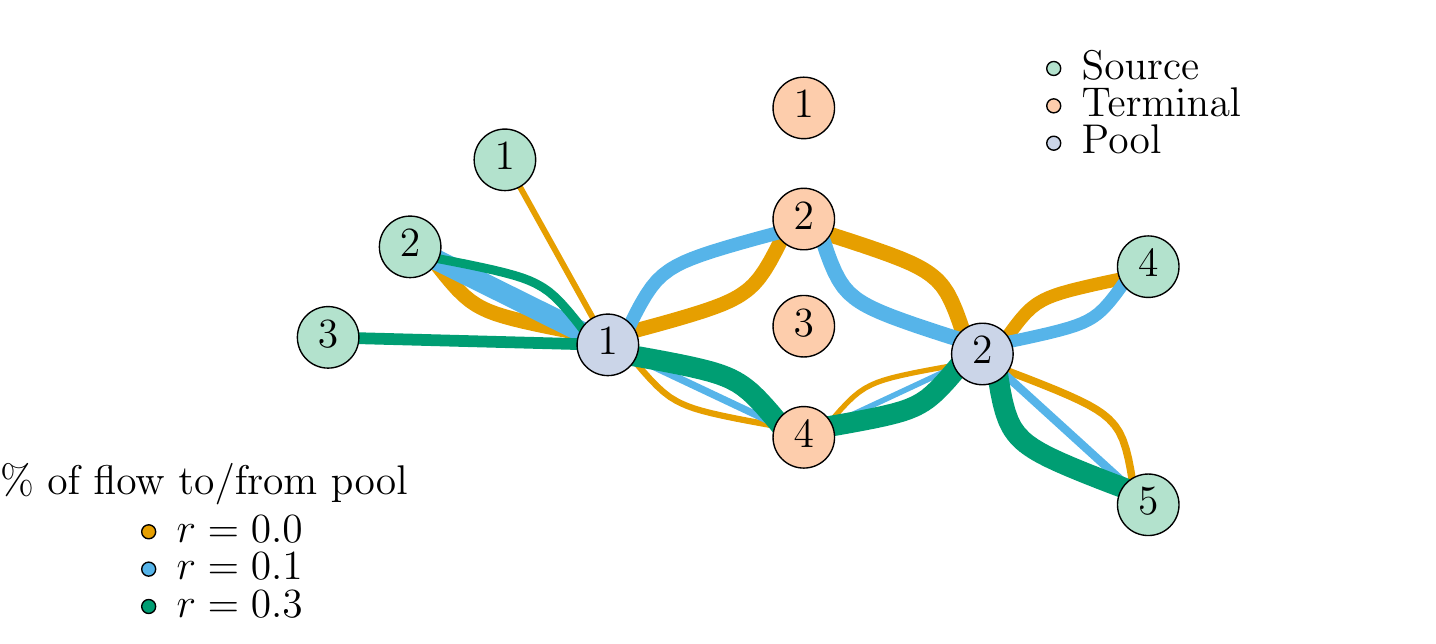}
    \caption{Active sources, pools, and terminals for three different uncertainy
        set sizes $r$ using instance Adhya~1 with the polyhedral uncertainty set.
        The width of the connecting arcs shows the percentage of flow to/from each pool.}
    \label{fig:adhya1}
    \vspace*{\floatsep}
    \includegraphics{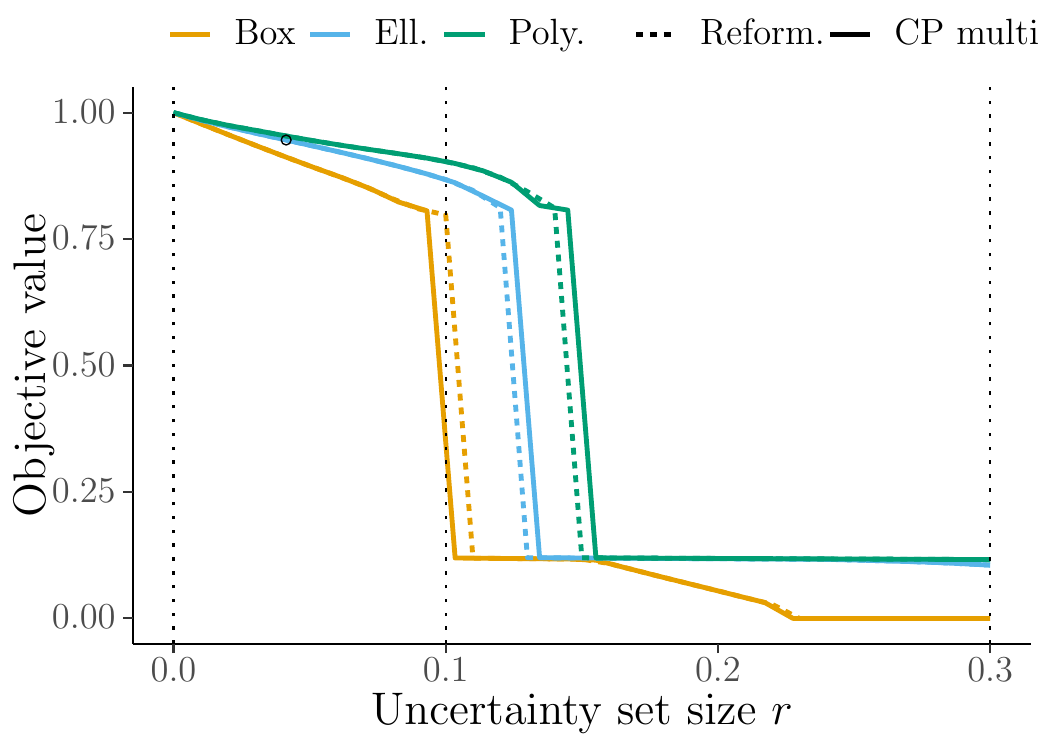}
    \caption{Objective value (relative to nominal case) for different
    uncertainty set types and sizes for instance Adhya~1. Black points indicate
instances which were not solved to global optimality within 1~hr.}
    \label{fig:adhya1-obj}
\end{figure}
Fig.~\ref{fig:adhya1} shows the optimal solution to instance
    Adhya~1 for three uncertainty set sizes $r = 0.0, 0.1, 0.3$ using the
    polyhedral uncertainty set.
    The width of the arrows between nodes
    indicates the fraction of flow from/to each pool.
In the nominal case, i.e., no uncertainty is considered ($r=0$),
Products~2 and~4 are produced using Sources~1, 2, 4, and~5.
As $r$ is increased, hedging against more uncertainty in the inlet
concentrations, at first the total
amount of Products~2 and~4 produced remains constant, but Source~1 is gradually
replaced by Source~2 which is more expensive but of better quality. When $r$
reaches $0.1$, Source~1 stops being useful altogether.
Past $r = 0.14$ the production of Product~2 ends, because its quality cannot be
guaranteed anymore for all parameter values in the uncertainty.
{\corr Stopping Product~2 production causes a sharp decrease in the expected profit from $446.2$
to $65.9$.}
Fig.~(\ref{fig:adhya1-obj}) shows the objective value of instance Adhya~1
(scaled by the nominal objective value) as a function of the uncertainty set
size $r$. The three scenarios shown in Fig.~(\ref{fig:adhya1}) and discussed
above are indicated by vertical dashed lines.
Results are shown for all three uncertainty set
geometries.
For a given uncertainty set geometry, as $r$ increases, the fraction of the
nominal objective value achieved always decreases, as would be expected.
Notice that, independent of the geometry, the profit first always decreases slowly
as cheap sources are substituted by higher quality sources and then
suddenly drops as the production of Product~2 becomes infeasible.
\setlength{\belowcaptionskip}{-4pt}
\begin{figure}[htb]
    \centering
    \includegraphics{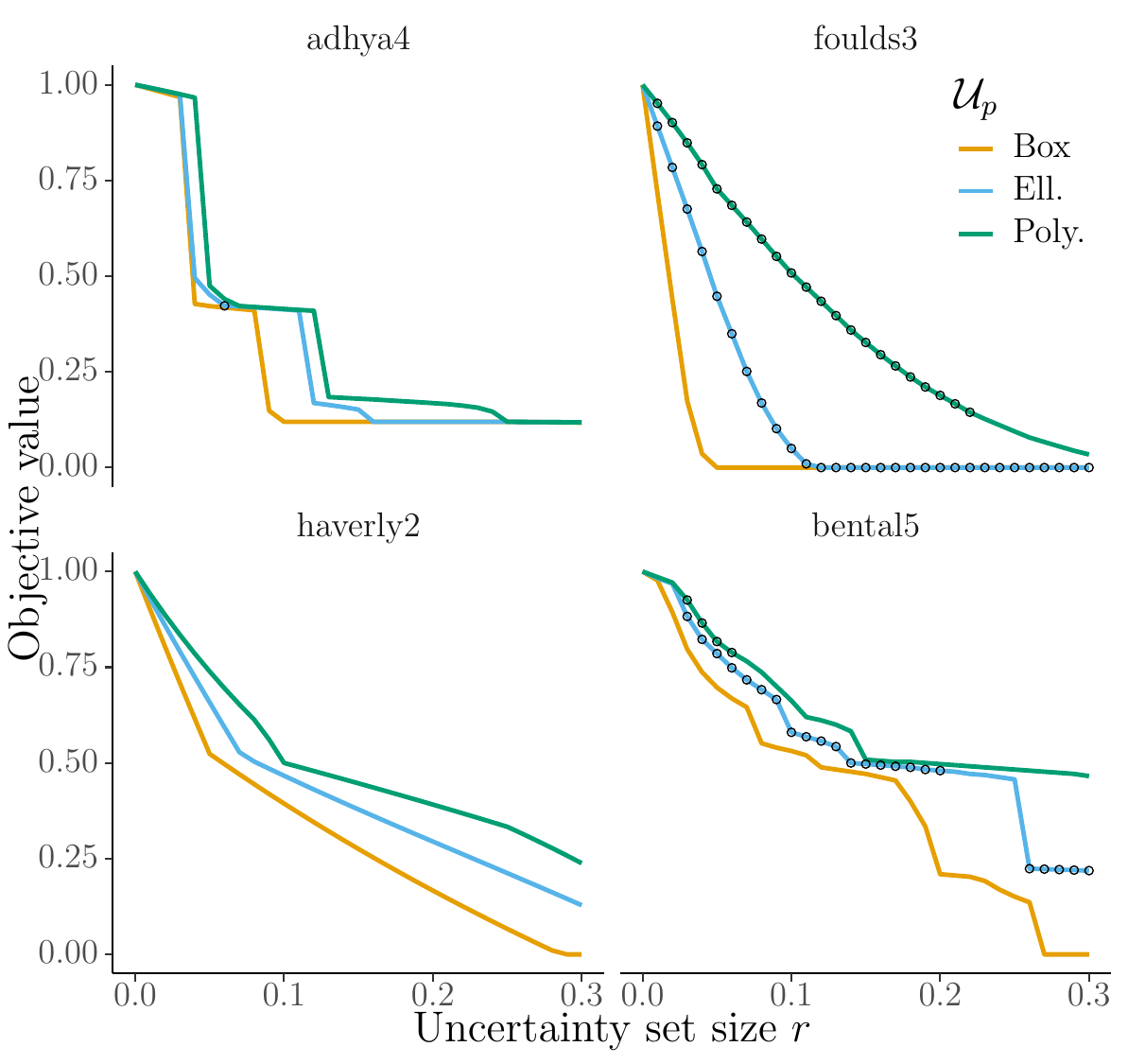}
    \caption{Objective value (relative to nominal case) for different
    uncertainty set types and sizes for four literature instances.
    Black points indicate instances which were not solved to global optimality within 1~hr.}
    \label{fig:obj-vs-r}
\end{figure}
Fig.~(\ref{fig:obj-vs-r}) shows similar trends for a number of other instances.

Figs.~(\ref{fig:adhya1-obj}) and~(\ref{fig:obj-vs-r}) show a clear ordering
between uncertainty sets geometries. By construction the
polyhedral uncertainty set is always smallest and the box uncertainty set
largest for a given $r$. Note that this does not necessarily mean that the
polyhedral set is always superior. While it is always less conservative, i.e.,
achieves a better worst case objective value, it may also be less robust.
Black points in Fig.~(\ref{fig:obj-vs-r}) indicate instances which could not be
solved to optimality within 1 hour.

\begin{figure}[htb]
    \centering
    \includegraphics{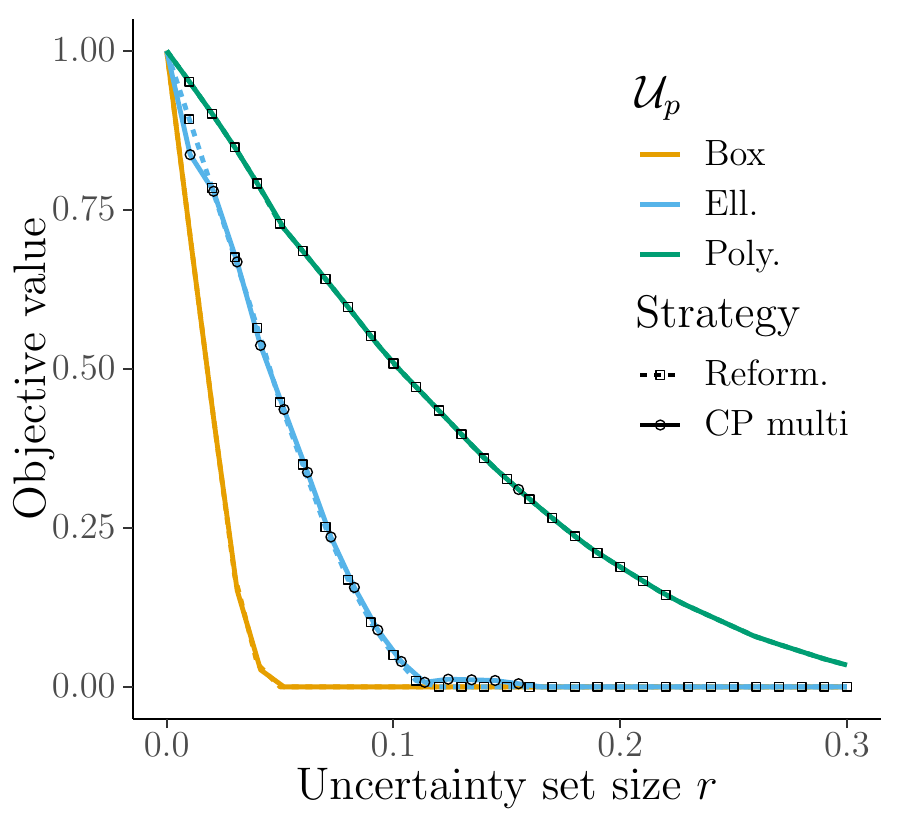}
    \caption{Objective value (relative to nominal case) as a function of
        uncertainty set size $r$ for instance Foulds~3. Results are shown for
        the reformulation approach and the multi-cut cutting plane approach.
        {\corr The black squares and circles indicate instances which were not solved
        to optimality within 1~hr by the reformulation and cutting plane
        strategy, respectively.} }
    \label{fig:cp}
\end{figure}
Fig.~(\ref{fig:cp}) shows a comparison for instance Foulds~3 between the
reformulation and the multi-cut version of the robust cutting plane approach.
As would be expected, the two approaches lead to identical objective values,
except for some small deviations with the ellipsoidal sets in regions where
neither approach converges to the optimal solution.
The cutting plane approach, however, converges for many more instances than the
reformulation approach. For this particular instance with the polyhedral
uncertainty set, it converges for all but
one value of $r$ while the reformulation approach only converges for very large
values of $r$.
A similar trend, albeit not as strong, can be seen for the ellipsoidal
uncertainty set.
\begin{table}[tb]
    \centering
    \small
    \begin{tabular}{l r r r r r r}
        \toprule
        Uncert. set & \multicolumn{3}{c}{\% instances solved} & \multicolumn{3}{c}{median time [s]}\\
                        & reform. & CP - multiple & CP - one & reform. & CP - multiple & CP - one \\
        \midrule
        Box          & 100 & 100 & 91 & 0.14 & 0.72 & 2.13 \\
        Ellipsoid    & 66  & 88  & 75 & 1.42 & 2.91 & 4.75 \\
        Polyhedron   & 83  & 96  & 79 & 0.15 & 1.13 & 1.77 \\
        \bottomrule
    \end{tabular}
    \caption{Percentage of instances solved to global optimality within 1~hr
    time limit and median time taken for different uncertainty set types and different robust
    optimization approaches. \corr{The median times were calculated for only those
    instances which could be solved to optimality within 1~hr.}}
    \label{tab:inst-solved}
\end{table}
\begin{figure}[htb]
    \centering
     \includegraphics{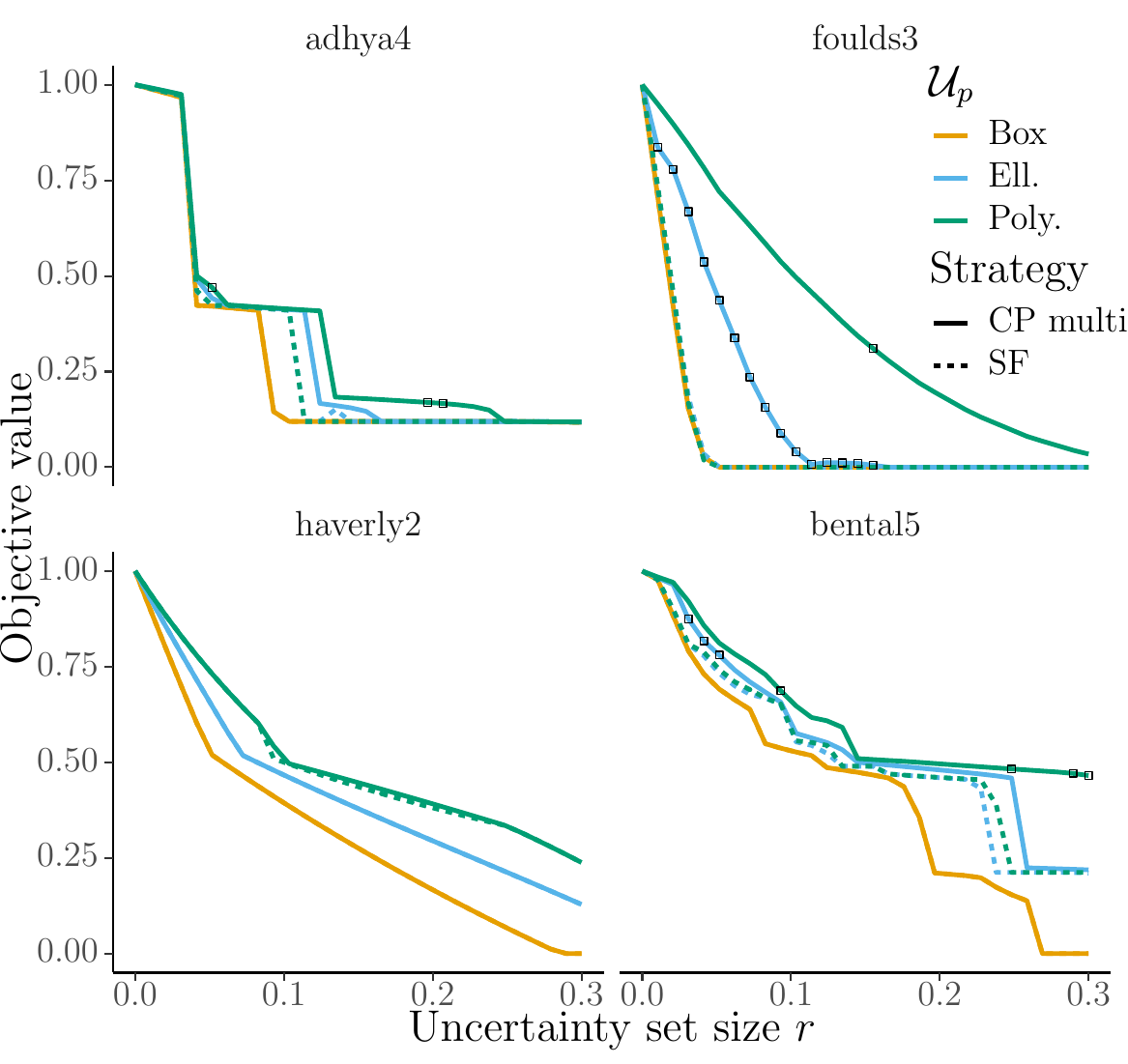}
    \caption{Objective value (relative to nominal case) as a function of
        uncertainty set size $r$ for the multi-cut cutting plane approach and
        the optimal safety factor approximation. Black points indicate
    instances which  were not solved to global optimality within 1~hr.}
    \label{fig:optimal-sf}
\end{figure}
Table~\ref{tab:inst-solved} shows this trend generally holds across instances.
The multi-cut cutting plane strategy solves a larger
percentage of instances than the reformulation approach for all uncertainty set
geometries. Note that the single-cut cutting plane approach, which adds only one
cut in every round, performs poorly in comparison to both alternatives.
The uncertainty set geometry also clearly has an effect
on the tractability of the robust problem. The box uncertainty set-constrained
problem is generally easy to solve. This is expected, especially for the
reformulation approach, because the problem difficulty is equivalent to
the nominal case. The reformulation approach performs particularly poorly for
the ellipsoidal set.
This is most likely due to the higher order non-linear terms in this formulation.
Table~\ref{tab:inst-solved} also shows the median time taken across
instances for those instances which could be solved to optimality.

{\corr While the cutting plane approach could
potentially add a large number of constraints, in practice it often adds a small number of constraints.
Across instances, the multi-cut strategy needs an
average of 13.3 iterations to converge. This leads to an average
final problem size of 335 constraints which is only slightly larger
than the average number of constraints in the reformulation approach:
328 for the polyhedral, 166 for the ellipsoidal, and 154 for the box
uncertainty set.} The cutting
plane approach does tends to take longer to solve than the reformulation approach.
This may, however, be at least partially because this approach solves a larger
fraction of the more difficult instances to optimality than the reformulation
approach.

While these results suggest that, at least for the pooling problem with the
uncertainty sets considered in this work, the cutting plane approach is
the preferable strategy for global robust optimization, the
reformulation approach has one important advantage: any feasible solution to
the robust reformulation is a valid feasible (albeit potentially non-optimal)
solution to the robust problem.
In contrast, intermediate solutions in the cutting plane algorithm are often
not robustly feasible.

Fig.~(\ref{fig:optimal-sf}) shows a comparison between the multi-cut robust
cutting plane approach and the safety factor approach outlined above.
For instance Haverly~2, the safety factor approach (with optimized safety
factors) is almost identical to the robust approach. This means that,
for this instance, even
the ellipsoidal and polyhedral uncertainty set on the inlet concentration
$\uncpar$ can be approximated by a box uncertainty set on $\PU$ and $\PL$.
However, e.g., for instances Adhya~4 and Bental~5, the safety-factor approach
is not equivalent anymore and leads to a conservative solution for both the
ellipsoidal and polyhedral set. For instance Foulds~3, the safety factor
approach is equivalent to the box uncertainty set. Overall, while the safety
factor is very effective in the sense that it does not increase the complexity
of the problem, it is generally a conservative approximation
to the robust problem.

{\corr
Fig.~(\ref{fig:correlation}) demonstrates correlation between sources.
The figure considers three different geographic layouts
for instance Adhya~1 with varying degrees of proximity between sources. The
for the ``far'' scenario are very similar to those of the ellipsoidal
uncertainty set without correlation. With decreasing distance between sources,
i.e., increasing correlation, the worst case profit generally decreases because
correlation between sources will make it more likely that, e.g., inlet concentrations
at multiple sources are all large simultaneously.
\begin{figure}[htb]
    \centering
    \includegraphics{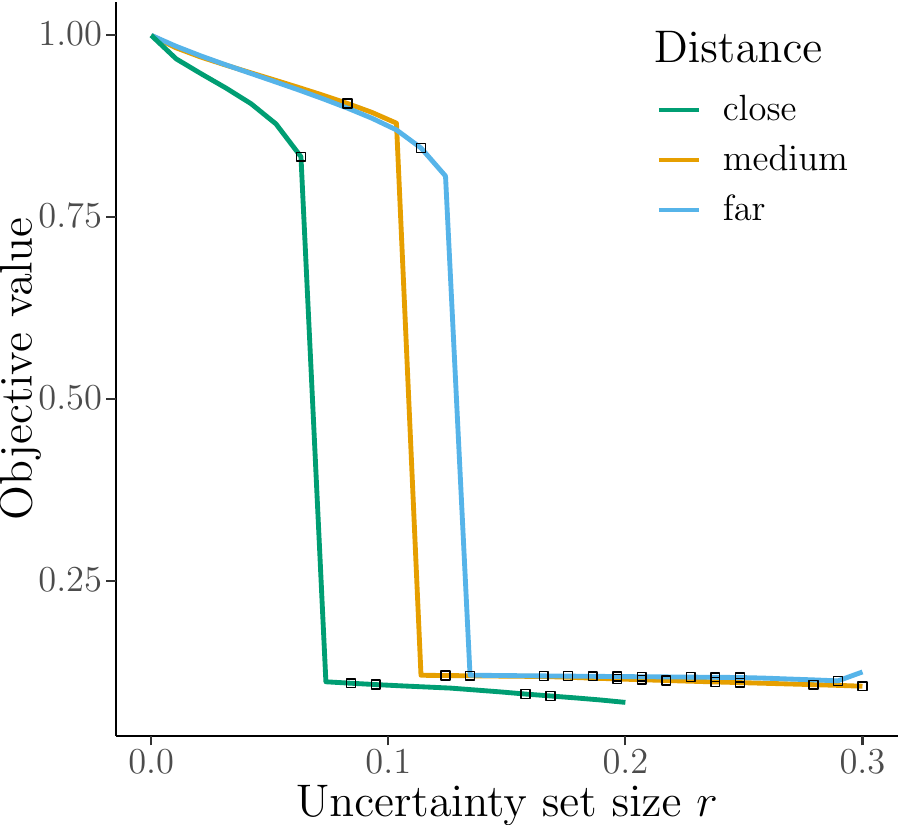}
    \caption{\corr Objective value (relative to nominal case) as a function of
        uncertainty set size $r$ for instance Adhya~1 with different geographic
        source locations using the multi-cut cutting plane approach.}
    \label{fig:correlation}
\end{figure}
}

\section{Conclusion}
The pooling problem is an interesting case study for robust optimization
because it is non-convex in the decision variables but concave in all
potentially uncertain parameters.
Applying robust reformulation or cutting plane approaches to this problem leads
to a global optimization problem with increased complexity compared to the nominal
case. The tractability of the resulting optimization problem with global optimization
solvers is highly dependent on the selected uncertainty set. While simple box
uncertainty sets hardly increase computational time, more advanced ellipsoidal
and polyhedral sets are more difficult to solve. For the pooling problem, the
cutting plane approach generally solves a larger percentage of instances to
optimality while the reformulation approach has the advantage that it produces
guaranteed feasible intermediate solutions. The relevance of optimization
under uncertainty is supported by the significant changes in the solution
observed when different degrees of uncertainty
are taken into account.

\section{Acknowledgements}
    This invited contribution is part of the I\&EC Research special issue for the 2019 Class of Influential
    Researchers.
    This work was funded by the Engineering \& Physical Sciences Research Council (EPSRC) Center for Doctoral Training in High Performance Embedded and Distributed Systems (EP/L016796/1), an EPSRC/Schlumberger CASE studentship to J.W. (EP/R511961/1, voucher 17000145), and an EPSRC Research Fellowship to R.M. (EP/P016871/1).

\bibliographystyle{compactnat}
\bibliography{ms}

\section{Nomenclature}

\begin{longtable}{ l l }
    \multicolumn{2}{l}{\textbf{\large Variables}}\\
    {$\vec{\lambda}$} & {vector of dual variables}\\
    {$\q$}          & {fraction of flow entering pool $l$ from source $i$}\\
    {$v_j$}          & {total flow to terminal $j$}\\
    {$\y$}          & {flow from pool $l$ to terminal $j$}\\
    {$\x$}        & {total flow from source $i$ to terminal $j$}\\
    {$\xvec$}       & {vector of decision variables}      \\[4pt]
    {$\xivec$}      & {vector of uncertain parameters} \\
    {$\z$}          & {direct flow from source $i$ to terminal $j$}\\
    \multicolumn{2}{l}{\textbf{\large Parameters}}\\
    {$\AL$}         & {lower availability limit feed $i$}\\
    {$\AU$}         & {upper availability limit feed $i$}\\
    {$b$}           & {right-hand-side constant}    \\
    {$c_i$}         & {per unit cost of material from source $i$}\\
    {$\uncpar$}     & {uncertain inlet concentration of component $k$ at source $i$}\\
    {$\nompar$}     & {nominal inlet concentration of component $k$ at source $i$}\\
    {$\devpar$}     & {maximum deviation of inlet concentration of component $k$ at source $i$}\\
    {$\delta$}      & {relative termination tolerance} \\
    {$d_j$}         & {per unit price of product at terminal $j$}\\
    {$\DL$}         & {minimum demand product $j$}\\
    {$\DU$}         & {maximum demand product $j$}\\
    {$\boldsymbol{\ell}_i$} & {location vector for source $i$}\\
    {$\PL$}         & {minimum concentration of component $k$ at terminal $j$}\\
    {$\PU$}         & {maximum concentration of component $k$ at terminal $j$}\\
    {$r$}           & {uncertainty set size parameter}\\
    {$S_l$}         & {maximum capacity pool $l$}\\
    {$s$}           & {safety factor}\\
    {$\sigma_{i,i'}$} & {covariance between input concentrations at sources $i$
    and $i'$}\\
    \multicolumn{2}{l}{\textbf{\large Sets}}\\
    {$\Lambda$}     & {feasible set of dual problem}\\
    {$T_X$}         & {contains $(i,l)$ if source $i$ is connected to pool $l$}\\
    {$T_Y$}         & {contains $(l,j)$ if pool $l$ is connected to terminal $j$}\\
    {$T_Z$}         & {contains $(i,j)$ if source $i$ is conntected to terminal $j$}\\
    {$\U$}          & {uncertainty set}\\
    {$\hat{\U}$}    & {set containing finite number of
                       uncertainty scenarios}\\
    {$\X$}          & {certain feasible set}\\
\end{longtable}

\end{document}

%% file: commands.tex
\renewcommand{\vec}{\boldsymbol}
\newcommand{\y}{y_{l,j}}
\newcommand{\q}{q_{i,l}}
\newcommand{\z}{z_{i,j}}
\newcommand{\x}{x_{i,j}}

\newcommand{\vj}{v_{j}}
\newcommand{\PU}{P^U_{j,k}}
\newcommand{\PL}{P^L_{j,k}}
\newcommand{\PUunc}{\tilde{P}^U_{j,k}}

\newcommand{\AU}{A^U_{i}}
\newcommand{\AL}{A^L_{i}}
\newcommand{\DU}{D^U_{j}}
\newcommand{\DL}{D^L_{j}}
\newcommand{\U}{\ensuremath{\mathcal{U}}}
\newcommand{\Up}{\ensuremath{\U_p}}
\newcommand{\uncpar}{\ensuremath{\tilde{C}_{i,k}}}
\newcommand{\nompar}{\ensuremath{{C}_{i,k}}      }
\newcommand{\devpar}{\ensuremath{\hat{C}_{i,k}}  }
\newcommand{\X}{\mathcal{X}}
\newcommand{\Ucp}{\hat{\U}}
\newcommand{\hp}{\lambda_{j,k}}
\newcommand{\T}{^\intercal}
\newcommand{\g}{g(\vec{x}, \vec{\xi})}
\newcommand{\f}{f(\vec{x})}
\newcommand{\xvec}{\vec{x}}
\newcommand{\xivec}{\vec{\xi}}
\newcommand{\smin}{s^{\text{min}}}
\DeclareMathOperator*{\argmax}{arg\,max}
\DeclareMathOperator*{\argmin}{arg\,min}
\newcommand{\R}{\mathbbm{R}}
\newcommand{\corr}{}

%% file: ms.bbl
\begin{thebibliography}{77}
\providecommand{\natexlab}[1]{#1}
\providecommand{\url}[1]{\texttt{#1}}
\expandafter\ifx\csname urlstyle\endcsname\relax
  \providecommand{\doi}[1]{doi: #1}\else
  \providecommand{\doi}{doi: \begingroup \urlstyle{rm}\Url}\fi

\bibitem[Janak and Floudas(2005)]{Janak2005}
Stacy~L. Janak and Christodoulos~A. Floudas.
\newblock {Advances in robust optimization approaches for scheduling under
  uncertainty}.
\newblock \emph{Computer Aided Chemical Engineering}, 20\penalty0 (C):\penalty0
  1051--1056, 2005.

\bibitem[Li and Ierapetritou(2008)]{Li2008}
Zukui Li and Marianthi~G. Ierapetritou.
\newblock {Robust Optimization for Process Scheduling Under Uncertainty}.
\newblock \emph{Industrial {\&} Engineering Chemistry Research}, 47\penalty0
  (12):\penalty0 4148--4157, 2008.

\bibitem[Ben-Tal et~al.(2009)Ben-Tal, {El Ghaoui}, and Nemirovski]{Bental2009}
Aharon Ben-Tal, Laurent {El Ghaoui}, and Arkadi Nemirovski.
\newblock \emph{{Robust Optimization}}.
\newblock Princeton University Press, 2009.

\bibitem[Gounaris et~al.(2013)Gounaris, Wiesemann, and Floudas]{Gounaris2013}
Chrysanthos~E Gounaris, Wolfram Wiesemann, and Christodoulos~A. Floudas.
\newblock {The Robust Capacitated Vehicle Routing Problem Under Demand
  Uncertainty}.
\newblock \emph{Operations Research}, 61\penalty0 (3):\penalty0 677--693, 2013.

\bibitem[Zhang et~al.(2015)Zhang, Grossmann, Heuberger, Sundaramoorthy, and
  Pinto]{Zhang2015}
Qi~Zhang, Ignacio~E. Grossmann, Clara~F. Heuberger, Arul Sundaramoorthy, and
  Jose~M. Pinto.
\newblock {Air separation with cryogenic energy storage: Optimal scheduling
  considering electric energy and reserve markets}.
\newblock \emph{AIChE Journal}, 61\penalty0 (5):\penalty0 1547--1558, 2015.

\bibitem[Vujanic et~al.(2016)Vujanic, Goulart, and Morari]{Vujanic2016}
Robin Vujanic, Paul Goulart, and Manfred Morari.
\newblock {Robust Optimization of Schedules Affected by Uncertain Events}.
\newblock \emph{Journal of Optimization Theory \& Applications}, 171\penalty0
  (3):\penalty0 1033--1054, 2016.

\bibitem[Ning and You(2017)]{Ning2017}
Chao Ning and Fengqi You.
\newblock {A data-driven multistage adaptive robust optimization framework for
  planning and scheduling under uncertainty}.
\newblock \emph{AIChE Journal}, 63\penalty0 (10):\penalty0 4343--4369, 2017.

\bibitem[Shang and You(2018)]{Shang2018}
Chao Shang and Fengqi You.
\newblock {Distributionally robust optimization for planning and scheduling
  under uncertainty}.
\newblock \emph{Computers \& Chemical Engineering}, 110:\penalty0 53--68, 2018.

\bibitem[Wiebe et~al.(2018)Wiebe, Cec{\'{i}}lio, and Misener]{Wiebe2018}
Johannes Wiebe, In{\^{e}}s Cec{\'{i}}lio, and Ruth Misener.
\newblock {Data-Driven Optimization of Processes with Degrading Equipment}.
\newblock \emph{Industrial {\&} Engineering Chemistry Research}, 57\penalty0
  (50):\penalty0 17177--17191, 2018.

\bibitem[Bertsimas et~al.(2010)Bertsimas, Brown, and Caramanis]{Bertsimas2010}
Dimitris Bertsimas, David B~Db Brown, and Constantine Caramanis.
\newblock {Theory and Applications of Robust Optimization}.
\newblock \emph{Operations Research}, page~50, 2010.

\bibitem[{Mohajerin Esfahani} and Kuhn(2017)]{MohajerinEsfahani2017a}
Peyman {Mohajerin Esfahani} and Daniel Kuhn.
\newblock \emph{{Data-driven distributionally robust optimization using the
  Wasserstein metric: performance guarantees and tractable reformulations}}.
\newblock Springer Berlin Heidelberg, 2017.

\bibitem[Bertsimas et~al.(2011)Bertsimas, Brown, and Caramanis]{Bertsimas2011}
Dimitris Bertsimas, David~B Brown, and Constantine Caramanis.
\newblock {Theory and Applications of Robust Optimization}.
\newblock \emph{SIAM Review}, 53\penalty0 (3):\penalty0 464--501, 2011.

\bibitem[Li et~al.(2011{\natexlab{a}})Li, Ding, and Floudas]{Li2011b}
Zukui Li, Ran Ding, and Christodoulos Floudas.
\newblock {A Comparative Theoretical and Computational Study on Robust
  Counterpart Optimization: I. Robust Linear Optimization and Robust Mixed
  Integer Linear Optimization}.
\newblock \emph{Industrial {\&} Engineering Chemistry Research}, 50\penalty0
  (18):\penalty0 10567--10603, 2011{\natexlab{a}}.

\bibitem[Li et~al.(2012{\natexlab{a}})Li, Tang, and Floudas]{Li2012b}
Zukui Li, Qiuhua Tang, and Christodoulos~A. Floudas.
\newblock {A Comparative Theoretical and Computational Study on Robust
  Counterpart Optimization: II. Probabilistic Guarantees on Constraint
  Satisfaction}.
\newblock \emph{Industrial {\&} Engineering Chemistry Research}, 51\penalty0
  (19):\penalty0 6769--6788, 2012{\natexlab{a}}.

\bibitem[Soyster(1973)]{Soyster}
A~L Soyster.
\newblock {Technical Note — Convex Programming with Set-Inclusive Constraints
  and Applications to Inexact Linear Programming}.
\newblock \emph{Operations Research}, 21\penalty0 (5):\penalty0 1154--1157,
  1973.

\bibitem[Ben-Tal and Nemirovski(1999)]{Bental1999}
A.~Ben-Tal and A.~Nemirovski.
\newblock {Robust solutions of uncertain linear programs}.
\newblock \emph{Operations Research Letters}, 25\penalty0 (1):\penalty0 1--13,
  1999.

\bibitem[Bertsimas and Sim(2004)]{Bertsimas2004}
Dimitris Bertsimas and Melvyn Sim.
\newblock {The Price of Robustness}.
\newblock \emph{Operations Research}, 52\penalty0 (1):\penalty0 35--53, 2004.

\bibitem[Ben-Tal and Nemirovski(1998)]{Bental1998}
A~Ben-Tal and A~Nemirovski.
\newblock {Robust Convex Optimization}.
\newblock \emph{Mathematics of Operations Research}, 23\penalty0 (4):\penalty0
  769--805, 1998.

\bibitem[Ben-Tal et~al.(2015)Ben-Tal, den Hertog, and Vial]{Bental2015}
Aharon Ben-Tal, Dick den Hertog, and Jean-Philippe Vial.
\newblock {Deriving robust counterparts of nonlinear uncertain inequalities}.
\newblock \emph{Mathematical Programming}, 149\penalty0 (1-2):\penalty0
  265--299, 2015.

\bibitem[Marandi et~al.(2017)Marandi, Ben-Tal, den Hertog, and
  Melenberg]{Marandi2017}
Ahmadreza Marandi, Aharon Ben-Tal, Dick den Hertog, and Bertrand Melenberg.
\newblock {Extending the Scope of Robust Quadratic Optimization}.
\newblock Technical report, Optimization Online, 2017.

\bibitem[Leyffer et~al.(2018)Leyffer, Menickelly, Munson, Vanaret, and
  Wild]{Leyffer2018}
Sven Leyffer, Matt Menickelly, Todd Munson, Charlie Vanaret, and Stefan~M Wild.
\newblock {Nonlinear Robust Optimization}.
\newblock 2018.
\newblock preprint ANL/MCS-P9040-0218, Argonne National Laboratory, Mathematics
  and Computer Science Division.

\bibitem[Misener and Floudas(2009)]{misener-floudas:2009}
R~Misener and CA~Floudas.
\newblock Advances for the pooling problem: Modeling, global optimization, and
  computational studies.
\newblock \emph{Appl. Comput. Math.}, 8\penalty0 (1):\penalty0 3 -- 22, 2009.

\bibitem[Misener et~al.(2011)Misener, Thompson, and Floudas]{Misener2011}
Ruth Misener, Jeffrey~P. Thompson, and Christodoulos~A. Floudas.
\newblock {APOGEE: Global optimization of standard, generalized, and extended
  pooling problems via linear and logarithmic partitioning schemes}.
\newblock \emph{Computers {\&} Chemical Engineering}, 35\penalty0 (5):\penalty0
  876--892, 2011.

\bibitem[Alfaki and Haugland(2013)]{Alfaki2013}
Mohammed Alfaki and Dag Haugland.
\newblock {Strong formulations for the pooling problem}.
\newblock \emph{Journal of Global Optimization}, 56\penalty0 (3):\penalty0
  897--916, 2013.

\bibitem[Floudas and Visweswaran(1990)]{Floudas1990}
C.A. Floudas and V.~Visweswaran.
\newblock {A global optimization algorithm (GOP) for certain classes of
  nonconvex NLPs -- I. Theory}.
\newblock \emph{Computers {\&} Chemical Engineering}, 14\penalty0
  (12):\penalty0 1397--1417, 1990.

\bibitem[Meyer and Floudas(2006)]{Meyer2006}
Clifford~A. Meyer and Christodoulos~A. Floudas.
\newblock {Global optimization of a combinatorially complex generalized pooling
  problem}.
\newblock \emph{AIChE Journal}, 52\penalty0 (3):\penalty0 1027--1037, 2006.

\bibitem[Quesada and Grossmann(1995)]{Quesada1995}
I.~Quesada and I.E. Grossmann.
\newblock {Global optimization of bilinear process networks with multicomponent
  flows}.
\newblock \emph{Computers {\&} Chemical Engineering}, 19\penalty0
  (12):\penalty0 1219--1242, 1995.

\bibitem[Tawarmalani and Sahinidis(2004)]{Tawarmalani2004}
Mohit Tawarmalani and Nikolaos~V. Sahinidis.
\newblock {Global optimization of mixed-integer nonlinear programs: A
  theoretical and computational study}.
\newblock \emph{Mathematical Programming}, 99\penalty0 (3):\penalty0 563--591,
  2004.

\bibitem[Wicaksono and Karimi(2008)]{Wicaksono2008}
Danan~Suryo Wicaksono and I.~A. Karimi.
\newblock {Piecewise MILP under- and overestimators for global optimization of
  bilinear programs}.
\newblock \emph{AIChE Journal}, 54\penalty0 (4):\penalty0 991--1008, 2008.

\bibitem[Gounaris et~al.(2009)Gounaris, Misener, and Floudas]{Gounaris2009}
Chrysanthos~E. Gounaris, Ruth Misener, and Christodoulos~A. Floudas.
\newblock {Computational Comparison of Piecewise−Linear Relaxations for
  Pooling Problems}.
\newblock \emph{Industrial {\&} Engineering Chemistry Research}, 48\penalty0
  (12):\penalty0 5742--5766, 2009.

\bibitem[Misener and Floudas(2010)]{Misener2010}
Ruth Misener and Christodoulos~A. Floudas.
\newblock {Global Optimization of Large-Scale Generalized Pooling Problems:
  Quadratically Constrained MINLP Models}.
\newblock \emph{Industrial {\&} Engineering Chemistry Research}, 49\penalty0
  (11):\penalty0 5424--5438, 2010.

\bibitem[Ceccon et~al.(2016)Ceccon, Kouyialis, and Misener]{Ceccon2016}
Francesco Ceccon, Georgia Kouyialis, and Ruth Misener.
\newblock {Using Functional Programming to Recognize Named Structure in an
  Optimization Problem: Application to Pooling}.
\newblock \emph{AIChE Journal}, 62\penalty0 (9):\penalty0 3085--3095, 2016.

\bibitem[Baltean-Lugojan and Misener(2018)]{Baltean2018}
Radu Baltean-Lugojan and Ruth Misener.
\newblock {Piecewise parametric structure in the pooling problem: from sparse
  strongly-polynomial solutions to NP-hardness}.
\newblock \emph{Journal of Global Optimization}, 71\penalty0 (4):\penalty0
  655--690, 2018.

\bibitem[Li et~al.(2012{\natexlab{b}})Li, Misener, and Floudas]{Li2012a}
Jie Li, Ruth Misener, and Christodoulos~A. Floudas.
\newblock {Scheduling of crude oil operations under demand uncertainty: A
  robust optimization framework coupled with global optimization}.
\newblock \emph{AIChE Journal}, 58\penalty0 (8):\penalty0 2373--2396,
  2012{\natexlab{b}}.

\bibitem[Li et~al.(2011{\natexlab{b}})Li, Armagan, Tomasgard, and
  Barton]{Li2011}
Xiang Li, Emre Armagan, Asgeir Tomasgard, and Paul~I. Barton.
\newblock {Stochastic pooling problem for natural gas production network design
  and operation under uncertainty}.
\newblock \emph{AIChE Journal}, 57\penalty0 (8):\penalty0 2120--2135,
  2011{\natexlab{b}}.

\bibitem[Li et~al.(2012{\natexlab{c}})Li, Tomasgard, and Barton]{Li2012}
Xiang Li, Asgeir Tomasgard, and Paul~I. Barton.
\newblock {Decomposition strategy for the stochastic pooling problem}.
\newblock \emph{Journal of Global Optimization}, 54\penalty0 (4):\penalty0
  765--790, 2012{\natexlab{c}}.

\bibitem[Li(2013)]{Li2013}
Xiang Li.
\newblock {Parallel nonconvex generalized Benders decomposition for natural gas
  production network planning under uncertainty}.
\newblock \emph{Computers \& Chemical Engineering}, 55:\penalty0 97--108, 2013.

\bibitem[Li and Barton(2015)]{Li2015}
Xiang Li and Paul~I Barton.
\newblock {Optimal design and operation of energy systems under uncertainty}.
\newblock \emph{Journal of Process Control}, 30:\penalty0 1--9, 2015.

\bibitem[Yang et~al.(2017)Yang, Vayanos, and Barton]{Yang2017}
Yu~Yang, Phebe Vayanos, and Paul~I Barton.
\newblock {Chance-Constrained Optimization for Refinery Blend Planning under
  Uncertainty}.
\newblock \emph{Industrial \& Engineering Chemistry Research}, 56\penalty0
  (42):\penalty0 12139--12150, 2017.

\bibitem[Kannan(2018)]{kannan-phd}
Rohit Kannan.
\newblock \emph{Algorithms, analysis and software for the global optimization
  of two-stage stochastic programs}.
\newblock PhD thesis, Massachusetts Institute of Technology, 2018.

\bibitem[Seong et~al.(2014)Seong, Chachuat, and Shah]{Seong2014}
Cheng Seong, Benoit Chachuat, and Nilay Shah.
\newblock {Fixed- flowrate total water network synthesis under uncertainty with
  risk management}.
\newblock \emph{Journal of Cleaner Production}, 77:\penalty0 79--93, 2014.

\bibitem[Wiebe et~al.(2019)Wiebe, Cec{\'{i}}lio, and Misener]{Wiebe2019}
Johannes Wiebe, In{\^{e}}s Cec{\'{i}}lio, and Ruth Misener.
\newblock {The robust pooling problem}.
\newblock In \emph{Proceedings of the 29th European Symposium on Computer Aided
  Process Engineering - ESCAPE 29}, 2019.
\newblock \textbf{accepted}.

\bibitem[Yuan et~al.(2017{\natexlab{a}})Yuan, Li, and Huang]{Yuan2016}
Yuan Yuan, Zukui Li, and Biao Huang.
\newblock {Robust optimization approximation for joint chance constrained
  optimization problem}.
\newblock \emph{Journal of Global Optimization}, 67\penalty0 (4):\penalty0
  805--827, 2017{\natexlab{a}}.

\bibitem[Polak(2012)]{Polak2012}
Elijah Polak.
\newblock \emph{Optimization: algorithms and consistent approximations}.
\newblock Springer Science \& Business Media, 2012.

\bibitem[Ben-Tal et~al.(2014)Ben-Tal, den Hertog, and Vial]{Ben-Tal2014}
Aharon Ben-Tal, Dick den Hertog, and Jean~Philippe Vial.
\newblock {Deriving robust counterparts of nonlinear uncertain inequalities}.
\newblock \emph{Mathematical Programming}, 149\penalty0 (1-2):\penalty0
  265--299, 2014.

\bibitem[Tsoukalas et~al.(2009)Tsoukalas, Rustem, and
  Pistikopoulos]{Tsoukalas2009}
Angelos Tsoukalas, Ber{\c{c}} Rustem, and Efstratios~N Pistikopoulos.
\newblock {A global optimization algorithm for generalized semi-infinite ,
  continuous minimax with coupled constraints and bi-level problems}.
\newblock \emph{Journal of Global Optimization}, 44:\penalty0 235--250, 2009.

\bibitem[Mitsos(2011)]{Mitsos2011}
Alexander Mitsos.
\newblock {Global optimization of semi-infinite programs via restriction of the
  right-hand side}.
\newblock \emph{Optimization}, 1934, 2011.

\bibitem[Stein(2012)]{Stein2012}
Oliver Stein.
\newblock {How to solve a semi-infinite optimization problem}.
\newblock \emph{European Journal of Operational Research}, 223\penalty0
  (2):\penalty0 312--320, 2012.

\bibitem[Diehl et~al.(2006)Diehl, Bock, and Kostina]{Diehl2006}
Moritz Diehl, Hans~Georg Bock, and Ekaterina Kostina.
\newblock {An approximation technique for robust nonlinear optimization}.
\newblock \emph{Mathematical Programming}, 107\penalty0 (1-2):\penalty0
  213--230, 2006.

\bibitem[Houska and Diehl(2013)]{Houska2013}
Boris Houska and Moritz Diehl.
\newblock {Nonlinear robust optimization via sequential convex bilevel
  programming}.
\newblock \emph{Mathematical Programming}, 142\penalty0 (1-2):\penalty0
  539--577, 2013.

\bibitem[Yuan et~al.(2017{\natexlab{b}})Yuan, Li, and Huang]{Yuan2017}
Yuan Yuan, Zukui Li, and Biao Huang.
\newblock {Nonlinear robust optimization for process design}.
\newblock \emph{AIChE Journal}, 64\penalty0 (2), 2017{\natexlab{b}}.

\bibitem[Ho-Nguyen and K{\i}l{\i}n{\c{c}}-Karzan(2018)]{Ho-nguyen2018}
Nam Ho-Nguyen and Fatma K{\i}l{\i}n{\c{c}}-Karzan.
\newblock {Exploiting problem structure in optimization under uncertainty via
  online convex optimization}.
\newblock \emph{Mathematical Programming}, 2018.

\bibitem[Stein and Still(2003)]{Stein2003}
Oliver Stein and Georg Still.
\newblock {Solving Semi-Infinite Optimization Problems with Interior Point
  Techniques}.
\newblock \emph{SIAM Journal on Control and Optimization}, 42\penalty0
  (3):\penalty0 769--788, 2003.

\bibitem[Diehl et~al.(2012)Diehl, Houska, Stein, and Steuermann]{Houska2012}
M~Diehl, B~Houska, O~Stein, and P~Steuermann.
\newblock {A lifting method for generalized semi-infinite programs based on
  lower level Wolfe duality}.
\newblock \emph{Computational Optimization and Applications}, 54\penalty0
  (1):\penalty0 189--210, 2012.

\bibitem[Harwood and Barton(2016)]{Harwood2016}
Stuart~M Harwood and Paul~I Barton.
\newblock {Lower level duality and the global solution of generalized
  semi-infinite programs}.
\newblock \emph{Optimization}, 65\penalty0 (6):\penalty0 1129--1149, 2016.

\bibitem[Dias and Ierapetritou(2016)]{Dias2016}
Lisia~S. Dias and Marianthi~G. Ierapetritou.
\newblock {Integration of scheduling and control under uncertainties: Review
  and challenges}.
\newblock \emph{Chemical Engineering Research \& Design}, 116:\penalty0
  98--113, 2016.

\bibitem[Castro et~al.(2018)Castro, Grossmann, and Zhang]{Castro2018}
Pedro~M. Castro, Ignacio~E. Grossmann, and Qi~Zhang.
\newblock {Expanding scope and computational challenges in process scheduling}.
\newblock \emph{Computers {\&} Chemical Engineering}, 114:\penalty0 14--42,
  2018.

\bibitem[Bertsimas et~al.(2016)Bertsimas, Dunning, and Lubin]{Bertsimas2015}
Dimitris Bertsimas, Iain Dunning, and Miles Lubin.
\newblock {Reformulation versus cutting-planes for robust optimization A
  computational study}.
\newblock \emph{Computational Management Science}, 13\penalty0 (2):\penalty0
  195--217, 2016.

\bibitem[Narraway and Perkins(1993)]{Narraway1993}
Lawrence~T Narraway and John~D Perkins.
\newblock {Selection of Process Control Structure Based on Linear Dynamic
  Economics}.
\newblock \emph{Industrial \& Engineering Chemistry Research}, 32\penalty0
  (11):\penalty0 2681--2692, 1993.

\bibitem[Bahri et~al.(2006)Bahri, Bandoni, and Romagnoli]{Bahri2006}
Parka~A Bahri, Jose~A Bandoni, and Jose~A Romagnoli.
\newblock {Effect of Disturbances in Optimizing Control : Steady-State
  Open-Loop Backoff Problem}.
\newblock \emph{AIChE Journal}, 42\penalty0 (4):\penalty0 983--994, 2006.

\bibitem[Koller et~al.(2018)Koller, Ricardez-Sandoval, and Biegler]{Koller2018}
Robert~W. Koller, Luis~A. Ricardez-Sandoval, and Lorenz~T. Biegler.
\newblock {Stochastic back-off algorithm for simultaneous design, control, and
  scheduling of multiproduct systems under uncertainty}.
\newblock \emph{AIChE Journal}, 64\penalty0 (7):\penalty0 2379--2389, 2018.

\bibitem[Rafiei and Ricardez-Sandoval(2018)]{Rafiei2018}
Mina Rafiei and Luis~A. Ricardez-Sandoval.
\newblock {Stochastic Back-Off Approach for Integration of Design and Control
  under Uncertainty}.
\newblock \emph{Industrial \& Engineering Chemistry Research}, 57\penalty0
  (12):\penalty0 4351--4365, 2018.

\bibitem[Haverly(1978)]{haverly}
C.~A. Haverly.
\newblock {Studies of the behavior of recursion for the pooling problem}.
\newblock \emph{ACM SIGMAP Bulletin}, \penalty0 (25):\penalty0 19--28, 1978.

\bibitem[Ben-Tal et~al.(1994)Ben-Tal, Eiger, and Gershovitz]{bental}
Aharon Ben-Tal, Gideon Eiger, and Vladimir Gershovitz.
\newblock {Global minimization by reducing the duality gap}.
\newblock \emph{Mathematical Programming}, 63\penalty0 (1-3):\penalty0
  193--212, 1994.

\bibitem[Gorissen et~al.(2015)Gorissen, Yanikoğlu, and den
  Hertog]{Gorissen2015}
Bram~L. Gorissen, Ihsan Yanikoğlu, and Dick den Hertog.
\newblock {A practical guide to robust optimization}.
\newblock \emph{Omega (United Kingdom)}, 53:\penalty0 124--137, 2015.

\bibitem[Lappas and Gounaris(2016)]{Lappas2016}
Nikolaos~H Lappas and Chrysanthos~E Gounaris.
\newblock {Multi-stage adjustable robust optimization for process scheduling
  under uncertainty}.
\newblock \emph{AIChE Journal}, 62\penalty0 (5):\penalty0 1646--1667, 2016.

\bibitem[Chiles(2011)]{Chiles2011}
Jean-Paul Chiles.
\newblock \emph{{Geostatistics: modeling spatial uncertainty}}.
\newblock Wiley-Blackwell, 2 edition, 2011.

\bibitem[Hettich and Kortanek(1993)]{Hettich1993}
R.; Hettich and K.~O. Kortanek.
\newblock {Semi-Infinite Programming: Theory, Methods, and Applications}.
\newblock \emph{SIAM Review}, 35\penalty0 (3):\penalty0 380--429, 1993.

\bibitem[Fischetti and Monaci(2012)]{Fischetti2012}
Matteo Fischetti and Michele Monaci.
\newblock {Cutting plane versus compact formulations for uncertain (integer)
  linear programs}.
\newblock \emph{Mathematical Programming Computation}, 4\penalty0 (3):\penalty0
  239--273, 2012.

\bibitem[{M{\'\i}nguez} and {Casero-Alonso}(2019)]{Minguez2019}
Roberto {M{\'\i}nguez} and V{\'\i}ctor {Casero-Alonso}.
\newblock {On the convergence of cutting-plane methods for robust optimization
  with ellipsoidal uncertainty sets}.
\newblock \emph{arXiv e-prints}, art. arXiv:1904.01244, 2019.

\bibitem[Dao-Thien and Massoud(1974)]{Dao1974}
My~Dao-Thien and M.~Massoud.
\newblock {On the Relation Between the Factor of Safety and Reliability}.
\newblock \emph{Journal of Engineering for Industry}, 96\penalty0 (3):\penalty0
  853--857, 1974.

\bibitem[Clausen et~al.(2006)Clausen, Ove, and Nilsson]{Clausen2006}
Jonas Clausen, Sven Ove, and Fred Nilsson.
\newblock {Generalizing the safety factor approach}.
\newblock \emph{Reliability Engineering {\&} System Safety}, 91\penalty0
  (8):\penalty0 964--973, 2006.

\bibitem[Hansson(2010)]{Hansson2010}
Sven~Ove Hansson.
\newblock {Promoting inherent safety}.
\newblock \emph{Process Safety and Environmental Protection}, 88\penalty0
  (3):\penalty0 168--172, 2010.

\bibitem[Adhya et~al.(1999)Adhya, Tawarmalani, and Sahinidis]{adhya}
Nilanjan Adhya, Mohit Tawarmalani, and Nikolaos~V. Sahinidis.
\newblock {A Lagrangian Approach to the Pooling Problem}.
\newblock \emph{Industrial {\&} Engineering Chemistry Research}, 38\penalty0
  (5):\penalty0 1956--1972, 1999.

\bibitem[Foulds et~al.(1992)Foulds, Haugland, and J{\"{o}}rnsten]{foulds}
L.~R. Foulds, D.~Haugland, and K.~J{\"{o}}rnsten.
\newblock {A bilinear approach to the pooling problem}.
\newblock \emph{Optimization}, 24\penalty0 (1-2):\penalty0 165--180, 1992.

\bibitem[Audet et~al.(2004)Audet, Brimberg, Hansen, Digabel, and
  Mladenovi{\'{c}}]{rt2}
Charles Audet, Jack Brimberg, Pierre Hansen, S{\'{e}}bastien~Le Digabel, and
  Nenad Mladenovi{\'{c}}.
\newblock {Pooling Problem: Alternate Formulations and Solution Methods}.
\newblock \emph{Management Science}, 50\penalty0 (6):\penalty0 761--776, 2004.

\bibitem[Misener and Floudas(2014)]{antigone}
Ruth Misener and Christodoulos~A. Floudas.
\newblock {ANTIGONE: Algorithms for coNTinuous / Integer Global Optimization of
  Nonlinear Equations}.
\newblock \emph{Journal of Global Optimization}, 59\penalty0 (2-3):\penalty0
  503--526, 2014.

\end{thebibliography}
